\documentclass[12pt,reqno]{amsart}
\usepackage{exscale}
\usepackage{relsize}
\usepackage{amsfonts}
\usepackage{a4wide}
\usepackage{mathrsfs,amsmath,makecell,amssymb}
\usepackage{makecell}
\usepackage{amssymb}
\usepackage{multicol}
\usepackage{color}
\allowdisplaybreaks
\numberwithin{equation}{section}
\newtheorem{theorem}{Theorem}[section]
\newtheorem{proposition}[theorem]{Proposition}

\newtheorem{lemma}[theorem]{Lemma}

\theoremstyle{definition}

\newtheorem{remark}[theorem]{Remark}

\newcommand{\R}{\mathbb{R}}
\newcommand{\dis}{\displaystyle}
\begin{document}

\title
 [Normalized solutions to $(2,q)$-Laplacian equation]
 {Normalized solutions to a class\\ of $(2,q)$-Laplacian equations}

\maketitle
\begin{center}

\author{Laura Baldelli}
\footnote{Email addresses: lbaldelli@impan.pl (L. Baldelli).}

\author{Tao Yang}
\footnote{Corresponding Author: Tao Yang. Email addresses: yangtao@zjnu.edu.cn (T. Yang).}
\end{center}

\begin{center}

\address {1 Institute of Mathematics, Polish Academy of Sciences,\\ ul. Sniadeckich 8, 00656, Warsaw, Poland}

\address {2 Department of Mathematics, Zhejiang Normal University, \\Jinhua, Zhejiang, 321004, P. R. China}

\end{center}
\maketitle

\begin{abstract}
This paper concerns the existence of normalized solutions to a class of $(2,q)$-Laplacian equations in all the possible cases according to the value of $p$ with respect to the critical exponent $2(1+2/N)$. In the $L^2$-subcritical case, we study a global minimization problem and obtain a ground state solution. While in the $L^2$-critical case, we prove several nonexistence results, extended also in the $L^q$-critical case. At last, we derive a ground state and infinitely many radial solutions in the $L^2$-supercritical case. Compared with the classical Schr\"{o}dinger equation, the $(2,q)$-Laplacian equation possesses a quasi-linear term, which brings in some new difficulties and requires a more subtle analysis technique. Moreover, the vector field $\vec{a}(\xi)=|\xi|^{q-2}\xi$ corresponding to the $q$-Laplacian is not strictly monotone when $q<2$, so we shall consider separately the case $q<2$ and the case $q>2$.


\end{abstract}

\maketitle

\section{Introduction and main result}

\setcounter{equation}{0}

This paper concerns the existence of solutions $(\lambda,u)\!\in\!  {\mathbb{R}}\!\times \!X$, with $X\!:=\!H^{1}({\mathbb{R}^N}) \!\cap\! D^{1,q}({\mathbb{R}^N})$, to the following $(2,q)$-Laplacian equation
\begin{equation}\label{eq1.1}
-\Delta u-\Delta_q u
   =\lambda u+ {| u |^{p - 2}}u \quad \text { in }\,\, \mathbb{R}^{N}
\end{equation}
under the constraint
\begin{equation}\label{eq1.2}
  \int_{{\mathbb{R}^N}} {{|u|}^2}dx=c^2,
\end{equation}
where $\Delta_q u\!=\!div({| \nabla u |^{q - 2}}\nabla u )$ is the $q$-Laplacian of $u$, $c\!>\!0$, $1\!<\!q\!<\!N$, $q\!\not=\! 2$, $2\!<\!p\!<\!\min\{2^{*}, q^*\}$ and $s^*\!:=\!\frac{sN}{N-s}$ is the critical Sobolev's exponent, for every $1<s<N$. In particular, we are seeking normalized solutions to \eqref{eq1.1}, since \eqref{eq1.2} imposes a normalization on its $L^2$-mass, which can be obtained by searching critical points of the following functional
\begin{equation} \label{energyF}
I(u)= \frac{1}{2} {\| {\nabla u} \|}_2^2  + \frac{1}{q}{\| \nabla u \|}_q^q- \frac{1}{p}{\|u \|}_p^p
\end{equation}
on the constraint
$$S_{c}: =\Big\{ u \in X:=H^{1}({\mathbb{R}^N}) \!\cap\! D^{1,q}({\mathbb{R}^N}): {\|u\|}_2^2=\int_{{\mathbb{R}^N}} {|u|}^2dx=c^2  \Big\}$$
with $\lambda$ appearing as Lagrange multipliers. This implies that $\lambda$ is part of the unknown.


The $(p, q)$-Laplacian equation \eqref{eq1.1} is closely related to the general reaction-diffusion system
\begin{align} \label{gpt1.4-1.5}
\partial_{t} u-\Delta_p u - \Delta_q  u = f(x, u),
\end{align}
which arises from chemical reactions, plasma physics, etc. Here, $u$ describes a concentration, the $(p, q)$-Laplacian term  in \eqref{gpt1.4-1.5} denote the diffusion as $div \big[ \big({| \nabla u |^{p - 2}}\!+\!{| \nabla u |^{q- 2}}\big)\nabla u\big]\!=\!\Delta_p u \!+\! \Delta_q  u$, $f(x, u)$ is the reaction related to sources and loss processes.
Another model related to the $(p, q)$-Laplacian operator concerns the study of solitary waves 
which is proposed in \cite{gDeK} modeling the elementary particles (cf. \cite{lByF,gDeK}). 

The stationary version of \eqref{gpt1.4-1.5} 
has been widely studied by many researchers. In \cite{CJLg}, C. J. He et al. proved the existence of a nontrivial solution under suitable assumptions on $f$. As a continuous work of \cite{CJLg}, they derived the regularity of weak solutions in \cite{HCl}. Later on, G. B. Li et al. \cite{LgZg} obtained multiple solutions with Sobolev's critical exponent.

Recently, L. Baldelli et al. \cite{lByF,lByF2} considered a critical problem of $(p,q)$-Laplacian type.
In particular, in \cite{lByF}, they obtained infinitely many weak solutions with negative energy in $D^{1,p}(\mathbb{R}^{N})\cap D^{1,q}(\mathbb{R}^{N})$ by using variational techniques and the concentration compactness principle by Lions. While in \cite{lByF2}, they studied the equation under a symmetric setting obtaining infinitely many solutions with positive energy.

While, the authors in \cite{VaDr} studied the effect of a potential term in a  $(p,q)$-Laplacian problem, namely
$$-\Delta_p u - \Delta_q  u +V(\varepsilon x)\big(|u|^{p-2}u+|u|^{q-2}u\big)=f(u)\quad \text { in }\,\, \mathbb{R}^{N},$$
where $\varepsilon\!>\!0$ is small, $1 \!<\!p\!<\!q\!<\!N$, $V\!\in\! C(\mathbb{R}^{N},\mathbb{R})$ satisfies the global Rabinowitz condition and $f\!\in \! C(\mathbb{R},\mathbb{R})$ is of subcritical growth. In detail, utilising the Ljusternik-Schnirelmann category theory, the authors derived the relation between the number of positive solutions and the topology of the set where $V$ attained its minimum for small $\varepsilon$. In \cite{NzgJ}, N. Zhang et al. obtained similar results to a class of $(p,q)$-Laplacian equations with sign-changing potential. Concerning $(p,q)$-Laplacian systems, the authors in \cite{RsZp} proved some existence and multiplicity results. For more information about $(p,q)$-Laplacian problems, please refer to \cite{VaDr,lByF,lByF2,LcYi,RsZp,HCl,CJLg,LgZg,NzgJ} and the references therein.

Taking the stationary version of \eqref{gpt1.4-1.5} with $p\!=\!2$, we obtain the $(2, q)$-Laplacian equation. The eigenvalue problem for a $(2, q)$-Laplacian equation is studied by V. Benci et al. in \cite{vBmA}, while, in \cite{nVrD,nVrD2,nVrD3}, N. S. Papageorgiou et al. proved several existence and multiplicity results to 
    $$-\Delta u - \Delta_q  u = f(x, u)~~~~\text{in}~~~~\Omega,\qquad u|_{\partial_ \Omega}=0$$
via variational methods and Morse theory if $q>2$.  

In this paper, motivated by the fact that physicists are often interested in normalized solutions, indeed prescribed mass appears in nonlinear optics and the theory of Bose-Einstein condensates, see \cite{fra, ma} and the reference therein,
we look for solutions of \eqref{eq1.1} in $X$ having a prescribed $L^2$-norm, as follows from \eqref{eq1.2}.
In literature, the existence of normalized solution $(\lambda, u)$
 to the following semilinear elliptic equation
\begin{equation}\label{n2}
-\Delta u=\lambda u+g(u), \qquad (\lambda, u)\in\mathbb R\times \mathbb R^N
\end{equation}
has been intensively studied recently. In the $L^2$-subcritical case, namely $p<2(1+\frac{2}{N})$, the functional on the constraint is coercive, so one can obtain the existence of a global minimizer by minimizing on the sphere, cfr. \cite{LilP, s82}, while in the other cases this method does not work. For instance, in the $L^2$-supercritical case, that is $p>2(1+\frac{2}{N})$, the functional on the sphere could not be bounded from below.
One of the main difficulties in dealing with normalized solutions as critical points of a functional constrained to a sphere consists in proving the Palais-Smale condition, as a compactness property.
Jeanjean in \cite{lJaE} overcomes these problems, in the $L^2$-supercritical case, using a mountain pass structure for an auxiliary functional proving the existence of at least one normalized solution of \eqref{n2}.
Later, the authors in \cite{TSdV} obtain infinitely many normalized solutions of \eqref{n2} using a new linking geometry for a stretched functional, by which we took inspiration for Theorem \ref{th1.4} below.
We also mention \cite{bm21} where, using a strong topological argument, the authors manage to avoid working with radial functions and Palais-Smale sequences, since the monotonicity of the ground state energy map is not required to deal with the lack of compactness.
Regarding the critical case, nonexistence results are reached in \cite{GZZ}, where the $p$-Laplacian operator is involved, together with a potential term, finally, we refer to \cite{Ww} where the critical Sobolev's exponent appears in the nonlinearity.
More results on normalized solutions for scalar equations can be found in \cite{BJL13, JLW, SJDE}, while for cooperative systems of coupled Schr\"odinger equations we refer to \cite{js22}.
See also \cite{LLY} where some results on normalized solutions to \eqref{n2} are well summarized.

Among papers dealing with more general equations, such as the nonlinear Choquard equation, we mention the paper by Bartsch et al. \cite{TlYl} where existence and multiplicity theorems are proved. See also \cite{Y20, LZ} for results in a fractional setting.
Moving to the quasilinear case, namely when the $p$-Laplacian operator is involved, very few is known in the local case, we cite \cite{Wlzk} where the supercritical case with a potential term is studied giving existence via a mountain pass argument. Concerning the nonlocal case, Hajaiej gives an important contribution to the field of normalized solutions to nonlinear Schr\"odinger equation with mixed fractional Laplacians of different types, see \cite{bgh, cgh, hp, lh}.

To the best of our knowledge, problem \eqref{eq1.1}-\eqref{eq1.2} has not been studied before. The goal of this paper is to start the analysis of the existence and nonexistence of solution to \eqref{eq1.1}-\eqref{eq1.2} in all the possible cases according to the value of $p$ with respect to the critical exponent $2(1+\frac{2}{N})$.
In particular, in the subcritical and supercritical case, we deal with ground state solutions, where $u$ is a \textbf{ground state} of \eqref{eq1.1} on $S_c$ if it is a solution to \eqref{eq1.1} having minimal energy among all solutions which belong to $S_c$, in other words
$$
(\left.I\right|_{S_{c}})'(u)=0 \quad \text { and } \quad I(u)=\inf \left\{I(w): (\left.I\right|_{S_{c}})'(w)=0 \text { and }  w\in S_{c}\right\}.
$$


In the $L^2$-subcritical case ($I|_{S_{c}}>-\infty$), we study a global minimization problem,
as stated below.

\begin{theorem}\label{th1.1}
Let $c>0$ and
\begin{equation}\label{h1.1a}
2<p<\biggl(1+\frac{2}{N}\biggl)\min\{2,q\},
\end{equation}
\begin{equation}\label{h1.1b}
\text{either}\qquad\frac{2N}{N+2}<q<2,\quad N\ge 2 \qquad \text{or}\qquad 2<q<N,\quad N\ge 3.
\end{equation}
Then $m(c)\!:=\!\inf _{u \in S_c} I(u)$ is achieved by some $u\in S_c$ with the following properties: $u$ is a nonnegative function in $\R^N$, is radially symmetric, solves \eqref{eq1.1} for some $\lambda_c$, and is a ground state of \eqref{eq1.1}-\eqref{eq1.2}. Moreover, we have
$$
-c^{\frac{2p(1-\delta_p)}{2-p\delta_p}} \lesssim m(c)<0,\qquad  -c^{\frac{2(p-2)}{2-p\delta_p}} \lesssim \lambda_c <0,
$$
hence $m(c) \to 0^{-}$ and $\lambda_c \to 0^{-}$ as $c \to 0^{+}$, where $\delta_p\!=\!\frac{N(p-2)}{2p}$.
\end{theorem}

Where, for any $a,b,c,d\in\mathbb R$, then $a \gtrsim b$ means $a \geq C b$ for some constant $C>0$. Similarly, $c \lesssim d$ means $c \leq C d$.

In the $L^2$-critical case, we obtain some nonexistence theorems.
\begin{theorem}\label{th1.2}
Let $N\!\geq\!2$ and $p\!=\!\bar{p}\!:=\!2(1\!+\!\frac{2}{N})$. Then there exists $ c_*\!>\!0$ such that   \\
(i) If $1\!<\!q\!<\!N$ and $0\!<\!c\!\leq\! c_*$, then $m(c)\!=\!0$ and there is no minimizer of $m(c)$;  \\
(ii) If $1\!<\!q\!<\!2$ and $ c\!>\!c_*$, then $m(c)\!=\!-\infty$ and there is no minimizer of $m(c)$.
\end{theorem}

In the following theorem, we analyze the $L^q$-critical case obtaining a partial result.

\begin{theorem}\label{th1.2bis}
Let $N>2$ and $p\!=\!\hat{p}\!:=q(1+\frac{2}{N})$. Then there exist $0<c_{**}<\hat{c}_{**}$ such that   \\
(i) If $1\!<\!q\!<\!N$ and $0\!<\!c\!\leq\! c_{**}$, then $m(c)\!=\!0$ and there is no minimizer of $m(c)$;  \\
(ii) If $2\!<\!q\!<\!N$ and $ c\!>\!\hat{c}_{**}$, then $m(c)\!=\!-\infty$ and there is no minimizer of $m(c)$.
\end{theorem}

Note that nonexistence in the interval $[c_{**},\hat{c}_{**}]$ remains an open problem since, as developed in Section \ref{cri}, asymptotic decays of particular externals seem to be missing in the literature, for results in this direction see \cite{bbo, cpy, Gdas}. 

In the $L^2$-supercritical case, since $I|_{S_{c}}$ is unbounded from below  ($I|_{S_{c}}=-\infty$), we consider a modified minimization problem
$$\sigma(c)\!:=\!\inf _{ u \in \mathcal{P}_{c} } I(u),$$
where
\begin{equation}  \label{eq2.4}
\mathcal{P}_{c}=\left\{u \in S_{c} : P(u)\!:=\! \|\nabla u\|_{2}^2+(1+\delta_{q})  \|\nabla u\|_{q}^q-\delta_{p}  {\|u\|}_p^p\!=\!0 \right\}.
\end{equation}
First, we prove the existence of a ground state (see Theorem \ref{th1.3}).

\begin{theorem}\label{th1.3}
Let $c>0$ and
\begin{equation}\label{h1.3a}
\biggl(1+\frac{2}{N}\biggl)\max\{2,q\}<p<\min\{2^*,q^*\},
\end{equation}
\begin{equation}\label{h1.3b}
\text{either}\,\,\,\,\frac{2N(N+2)}{N^2+2N+4}<q<2, \,\, N\ge 2 \quad \text{or}\quad 2<q<\min\biggl\{N, \frac{2N^2}{N^2-4}\biggr\}, \,\, N\ge 3
\end{equation}
Then $ \sigma(c)\!=\!\inf _{u \in \mathcal{P}_{c}} I(u)$ is attained by some $u\in \mathcal{P}_c$ with the following properties: $u$ is a nonnegative function in $\R^N$, is radially symmetric, solves \eqref{eq1.1} for some $\lambda_c<0$, and is a ground state of \eqref{eq1.1}-\eqref{eq1.2}. In addition, we have
$$\sigma(c) \gtrsim c^{-\frac{2p(1-\delta_{p})}{p\delta_{p}-2}}+c^{-\frac{qp(1-\nu_{p,q})}{p\nu_{p,q}-q}},\qquad \lambda_c \lesssim -\biggl(c^{-\frac{2(p-2)}{p\delta_p-2}}+c^{ -\frac{2q(p-2)}{Np-Nq-2q} }\biggr), $$
and hence $\sigma(c) \to +\infty$ and $\lambda_c \to -\infty$ as $c \to 0^{+}$, where $\nu_{p,q}\!=\!\frac{Nq(p-2)}{p[Nq-2(N-q)]}$.
 \end{theorem}

Next, we obtain infinitely many solutions in the radial space $X_{r}\!:=\!\{u(x) \!\in\! X: u(x)\!=\!u(|x|)\}$.

\begin{theorem}\label{th1.4}
Let $c>0$ and assume \eqref{h1.3a}, \eqref{h1.3b}. Then there exists $n_0\in \mathbb{N}^+$ such that, for any fixed $c>0$ and $n\geq n_0$, $\eqref{eq1.1}$ possesses a sequence of couples of weak solutions $\{({v^{(n)}},{\lambda^{(n)}})\}\subseteq X_r \times {\R^- }$ with $\left\| {{v^{(n)}}} \right\|_2=c$. Moreover, it holds that
$$
\|v^{(n)}\|_{X_r}\to +\infty,\qquad I(v^{(n)})\to +\infty~~~~\mbox{as}~~~~n \to +\infty.
$$
\end{theorem}

\begin{remark}\label{re1.3} As far as we know, this is the first result concerning the existence and multiplicity of normalized solutions to \eqref{eq1.1}. Therefore, we provide a new perspective to the study of the $(2, q)$-Laplacian equation. Our main results cover the $L^2$-subcritical case, $L^2$-critical case and $L^2$-supercritical case. Since ${\| {\nabla u} \|}_2^2$ and ${\| \nabla u \|}_q^q$ behave differently after scaling, we need to deal with the case $q>2$ and $q<2$ respectively. When $q>2$, the vector field $\vec{a}(\xi)=|\xi|^{q-2}\xi$ corresponding to the $q$-Laplacian is strictly monotone, that is to say $\vec{a}(\cdot)$ satisfies  $$\big(\vec{a}(\xi)-\vec{a}(\eta)\big)\big(\xi-\eta\big) \gtrsim |\xi-\eta|^q~~~~\text{for}~~~~\xi,\eta\in\R^N,$$
this property is vital in compactness analysis. However, since $\vec{a}(\xi)$ is not strictly monotone if $q<2$, we prove a useful inequality to tackle it.
\end{remark}

The $(2, q)$-Laplacian equation \eqref{eq1.1} involves a quasi-linear term, so the study of \eqref{eq1.1} is more difficult than that of the classical Schr\"{o}dinger equation 
$$- \Delta u=\lambda u+ {| u |^{p - 2}}u  \text { in } \mathbb{R}^{N}.$$ 
In fact, the workspace in studying \eqref{eq1.1} changes from $H^{1}({\mathbb{R}^N})$ to $H^{1}({\mathbb{R}^N}) \!\cap\! D^{1,q}({\mathbb{R}^N})$, the quasi-linear term cannot be precisely controlled and there exists a competition between $- \Delta$ and $- \Delta_q$ when $q\not=2$. Compared with the study of \eqref{eq1.1} when $\lambda$ is fixed, our case is more challenging because $\lambda$ is part of the unknown and the strong $L^2({\mathbb{R}^N})$ convergence of the selected Palais-Smale sequence in $X$ is hard to prove. \\

Now, we give the outline of the proofs to Theorems \ref{th1.1}-\ref{th1.4}. To prove Theorem \ref{th1.1}, we show that the global minimization problem $m(c)\!:=\!\inf _{u \in S_c} I(u)$ is attained. After have taken $\{v_{n}\} \subset S_c$, a minimizing sequence for $m(c)$, using Ekeland’s variational principle, we derive a new minimizing sequence $\{u_{n}\} \subset S_c$ that is also a Palais-Smale sequence for $I|_{S_c}$, such that
$$\left\|u_{n}-v_{n}\right\|_{X} \to 0,\quad \big(I|_{S_c}\big)'\left(u_{n}\right) \to 0\qquad \text{as}~~~~n\to+\infty.$$
Based on $\big(I|_{S_c}\big)'\left(u_{n}\right) \to 0$,  we can prove that there exists $u_{c}\in S_c$ such that $\nabla u_{n} \rightarrow \nabla u_c $ a.e. on ${\R}^N$, which is crucial in obtaining a Br\'{e}zis-Lieb splitting lemma. After excluding the vanishing of the weak limit of $\{u_{n}\}$, we utilize the strict subadditivity inequality $m(c)< m(c_1)+m(\sqrt{c^2-c^2_1})$, with $0<c_1<c$, to rule out the dichotomy of $\{u_{n}\}$. To this end, we take into account the symmetric decreasing rearrangement to conclude and prove the required inequality for $m(c)$ and $\lambda_c$. So, $m(c)$ has a minimizer in $S_c$. In Theorems \ref{th1.2}-\ref{th1.2bis}, we select some explicit test functions to detect when $I|_{S_{c}}>-\infty$ and $I|_{S_{c}}=-\infty$ and we analyse their properties.

In Theorems \ref{th1.3}-\ref{th1.4}, we are faced with a more complicated case where $I|_{S_{c}}=-\infty$. In this situation, it is a highly nontrivial issue to obtain bounded Palais-Smale sequences for $I|_{S_{c}}$. We tackle this difficulty by using a natural constraint approach in Theorem \ref{th1.3}, while in Theorem \ref{th1.4} we adopt a stretched functional. In particular, to prove Theorem \ref{th1.3}, we show that the modified minimization problem $ \sigma(c)\!:=\!\inf _{ u \in \mathcal{P}_{c} } I(u)$ possesses a  minimizer. The set $\mathcal{P}_{c}$ involves two constraints, one is a mass constraint and the other one is natural Pohozaev constraint $P(u)\!=\!0$ (see \eqref{eq2.4}). In view of this, we can check that $I|_{\mathcal{P}_{c}}>-\infty$ and hence $\sigma(c)$ makes sense. Next, we construct different scaling functions to show that $\sigma(c)$ is strictly decreasing with respect to $c$, which is a key ingredient in searching minimizers for $\sigma(c)$. This step is very technical due to the competition between $- \Delta$ and $- \Delta_q$. Then, combining the symmetric decreasing rearrangement technique with the monotonicity of $\sigma(c)$, we obtain a minimizer to $\sigma(c)$.

In proving Theorem \ref{th1.4}, the Kranoselski genus theory \cite{MsUe} is not applicable since $I|_{S_{c}}=-\infty$. We turn to use a linking argument of $I|_{S_{c}}$, which generates infinitely many radial solutions to \eqref{eq1.1}-\eqref{eq1.2}. To begin with, we prove a key intersection lemma (see Lemma \ref{6.1.inter}). Next, we adopt a stretched functional to create bounded Palais-Smale sequences for $I|_{S_{c}}$. Finally, we end the proof with a min-max procedure.   \\

The paper is organized as follows. In Section \ref{prel}, we give some preliminary results used in the proof of our theorems. Concerning the $L^2$-subcritical case, Section \ref{sub} encloses the proof Theorem \ref{th1.1} passing by intermediate Lemmas, while in Section \ref{cri} we prove nonexistence in the $L^2$ and $L^q$-critical cases, in terms of Theorems \ref{th1.2}-\ref{th1.2bis}. Finally, the $L^2$-supercritical case is described in Section \ref{supc} which is divided into two subsections: Subsection \ref{subsec5.1} is dedicated to proving Theorem \ref{th1.3} that gives the existence of a ground state solution, while Subsection \ref{subsec5.2} contains the proof of Theorem \ref{subsec5.2} devoted in the existence of infinitely many radial solutions.\\


\noindent \textbf{Notations:}~~~~For $1 \!\le\! p \!<\!\infty $ and $u\!\in\!{L^p}({\R^N})$, we denote ${\left\| u \right\|_p}\!:=\! {({\int_{{\R^N}} {\left| u \right|} ^p}dx)^{\frac{1}{p}}}$. The Hilbert space $H^{1}(\mathbb{R}^{N})$ is defined as $H^{1}(\mathbb{R}^{N}) :=\{u \in L^{2}(\mathbb{R}^{N}) :\nabla u \in L^{2}(\mathbb{R}^{N})\}$ with inner product $(u,v): = \int_{{\R^N}} (\nabla u\nabla v + uv)dx$  and norm ${\left\| u \right\|} := (\left\| {\nabla u} \right\|_2^2 + \left\| u \right\|_2^2)^{\frac{1}{2}}$. Similarly, $D^{1,q}(\mathbb{R}^{N})$ reads $D^{1,q}(\mathbb{R}^{N}) :=\{u \in L^{q^*}(\mathbb{R}^{N}) :\nabla u \in L^{q}(\mathbb{R}^{N})\}$ with the semi-norm $\|u\|_{D^{1,q}(\mathbb{R}^{N})}\!=\!\left\| {\nabla u} \right\|_q$. Recalling  $X=H^{1}({\mathbb{R}^N}) \!\cap\! D^{1,q}({\mathbb{R}^N})$ endowed with the norm $\|u\|_X=\|u\|+\|u\|_{D^{1,q}(\mathbb{R}^{N})}$, then $X^{-1}$ is the dual space of $X$. We use $``\rightarrow"$ and $``\rightharpoonup"$ to denote the strong and weak convergence in the related function spaces respectively. $C$ and $C_{i}$ will be positive constants. $\langle\cdot,\cdot\rangle$ denote the dual pair for any Banach space and its dual space.
Finally, $o_{n}(1)$ and $O_{n}(1)$ mean that $|o_{n}(1)|\to 0$ and $|O_{n}(1)|\leq C$ as $n\to+\infty$, respectively. 


\section{Preliminaries}\label{prel}

\setcounter{equation}{0}
In this section, we introduce various preliminary results.
\begin{lemma} (Gagliardo-Nirenberg inequality, \cite{Wein}) \label{lem2.1}
Let $p\!\in\!(2,2^*)$ and $\delta_p\!=\!\frac{N(p-2)}{2p}$. Then there exists a constant $\mathcal{C}_{N,p}\!=\!\Big( \frac{p}{2 \|W_p\|^{p-2}_{2}} \Big)^{\frac{1}{p}}\!>\!0$ such that
\begin{equation} \label{equ2.2}
 \|u\|_{p} \leq \mathcal{C}_{N,p} \left\|\nabla u\right\|_{2}^{\delta_p} \left\|u\right\|_{2}^{(1-\delta_p)}, \qquad \forall u \in {H}^{1}(\mathbb{R}^{N}),
\end{equation}
where $W_p$ is the unique positive radial solution of
$ -\Delta W\!+\!(\frac{1}{\delta_p}-\!1)W \!=\!\frac{2}{p\delta_p}|W|^{p-2}W$.
\end{lemma}

\begin{lemma} ($L^q$-Gagliardo-Nirenberg inequality, \cite[Theorem 2.1]{MaEh}) \label{lem2.2}
Let $q\!\in\!(\frac{2N}{N+2},N)$, $p\!\in\!(2,q^*)$ and $\nu_{p,q}$ as in Theorem \ref{th1.3}. Then there exists a constant $\mathcal{K}_{N,p}\!>\!0$ such that
\begin{equation} \label{pq-equ2.2}
 \|u\|_{p} \leq \mathcal{K}_{N,p} \left\|\nabla u\right\|_{q}^{\nu_{p,q}} \left\|u\right\|_{2}^{(1-\nu_{p,q})}, \qquad \forall u \in {D}^{1,q}(\mathbb{R}^{N})\cap{L}^{2}(\mathbb{R}^{N}),
\end{equation}
where 
$$\mathcal{K}_{N,p}=\biggl[\frac{K}{\frac{1}{q}\|DW_{p,q}\|_q^q+\frac{1}{2}\|W_{p,q}\|_2^2}\biggr],$$
$$K=(Nq+pq-2N)\cdot\biggr[\frac{[2(Nq-p(N-q))]^{p(N-q)-Nq}}{[qN(p-2)]^{N(p-2)}}\biggr]^{1/[Nq+pq-2N]},$$
and $W_{p,q}$ is the unique nonnegative radial solution of the following equation
$$-\Delta_q W+W=\zeta|W|^{p-2}W,$$
where $\zeta=\|DW\|_q^q+\|W\|_2^2$ is the Lagrangian multiplier.

\end{lemma}

\begin{lemma}  \label{lem2.7}
Let $N\!\ge\!2$, $q \!\in\!(1,N)$, $2\!<\!p\!<\!2^*$ and $\lambda\in\R$. If $u \!\in\! X $ is a weak solution of \eqref{eq1.1},
then the Pohozaev identity
$P(u)\!=\! \|\nabla u\|_{2}^2+(1+\delta_{q})  \|\nabla u\|_{q}^q-\delta_{p}  {\|u\|}_p^p\!=\!0$ holds. 
\end{lemma}
\begin{proof}
For simplicity, we rewrite ${u}_{x_i}=\frac{\partial u}{\partial x_i } $ and ${u}_{x_ix_j}=\frac{\partial}{\partial x_j}\Big(\frac{\partial u}{\partial x_i } \Big)$ and omit the integral symbol $dx$. We first suppose $u \!\in\! C^2(\mathbb{R}^{N})$. Integrating by parts, we derive
\begin{align*} 
\int_{\mathbb{R}^{N}}&(\Delta u)(x\cdot \nabla  {u}) \!=\!\sum_{i,j=1}^{N}\int_{\mathbb{R}^{N}}u_{x_ix_i}x_j {u}_{x_j}
\!=\!-\sum_{i,j=1}^{N}\int_{\mathbb{R}^{N}}u_{x_i}[\delta_{i,j}
 {u}_{x_j}+x_j {u}_{x_ix_j}]   \nonumber  \\
&\!=\!-\int_{\mathbb{R}^{N}} {| \nabla u|}^2-\sum_{i,j=1}^{N}\int_{\mathbb{R}^{N}}u_{x_i}x_j {u}_{x_ix_j}
\!=\!-\int_{\mathbb{R}^{N}} {| \nabla u|}^2-\sum_{i,j=1}^{N}\int_{\mathbb{R}^{N}}  \Big( \frac{|u_{x_i}|^2}{2}\Big)_{x_j}x_j
   \nonumber \\
&\!=\!-\int_{\mathbb{R}^{N}} {| \nabla u|}^2+\sum_{j=1}^{N}\int_{\mathbb{R}^{N}} \frac{|\nabla u |^2}{2}
\!=\! \frac{N-2}{2}\int_{\mathbb{R}^{N}} {| \nabla u|}^2
\end{align*}
and
\begin{align*} 
\int_{\mathbb{R}^{N}}(\Delta_q u)(x\cdot \nabla  {u}) &\!=\!\sum_{j=1}^{N}\int_{\mathbb{R}^{N}} div({| \nabla u |^{q - 2}}\nabla u) x_j {u}_{x_j}
\!=\!\sum_{i,j=1}^{N}\int_{\mathbb{R}^{N}}  ({| \nabla u |^{q - 2}}u_{x_i})_{x_i} x_j {u}_{x_j}   \nonumber  \\
&\!=\!-\sum_{i,j=1}^{N}\int_{\mathbb{R}^{N}}  ({| \nabla u |^{q - 2}}u_{x_i}) ( \delta_{i,j} {u}_{x_j}+x_j {u}_{x_ix_j})  \nonumber  \\
&\!=\!-\sum_{i=1}^{N}\int_{\mathbb{R}^{N}}  {| \nabla u |^{q - 2}}u_{x_i} {u}_{x_i}-\sum_{i,j=1}^{N}\int_{\mathbb{R}^{N}}  {| \nabla u |^{q - 2}}u_{x_i} x_j {u}_{x_ix_j}
   \nonumber \\
&\!=\!-\int_{\mathbb{R}^{N}} {| \nabla u|}^q-\sum_{j=1}^{N}\int_{\mathbb{R}^{N}}\Big( \frac{|\nabla u |^q}{q}\Big)_{x_j} x_j
   \!=\! \frac{N-q}{q}\int_{\mathbb{R}^{N}} {| \nabla u|}^q,
\end{align*}
where $\delta_{i,j}=1$ for $i=j$ and $\delta_{i,j}=0$ for $i\not=j$.
In the same way, for any $p\geq2$, we have
\begin{align*}
\int_{\mathbb{R}^{N}} {{|u|}^{p-2}u}(x\cdot \nabla {u})
=\sum_{j=1}^{N}\int_{\mathbb{R}^{N}} {{|u|}^{p-2}u} x_j {u}_{x_j}
=\sum_{j=1}^{N}\int_{\mathbb{R}^{N}} \Big( \frac{| u |^p}{p}\Big)_{x_j} x_j
=-\frac{N}{p}\int_{\mathbb{R}^{N}} {{|u|}^p}.
\end{align*}
Therefore, we have
\begin{align}\label{p3.4}
\frac{2-N}{2}{\| {\nabla u} \|}_2^2+\frac{q-N}{q}{\| {\nabla u} \|}_q^q\!=\!\frac{-N\lambda}{2} {{{\| u \|}_2^2}}-\frac{N}{p}{ {{\| u \|}_p^p}}.
\end{align}
To eliminate $\lambda$, we multiply \eqref{eq1.1} by ${u}$ and get
\begin{align}\label{p3.5}
{\| {\nabla u} \|}_2^2+{\| {\nabla u} \|}_q^q\!=\!\lambda {{{\| u \|}_2^2}}+{ {{\| u \|}_p^p}}.
\end{align}
Combining \eqref{p3.4}-\eqref{p3.5}, we have the Pohozaev identity $\|\nabla u\|_{2}^2+(1+\delta_{q})  \|\nabla u\|_{q}^q-\delta_{p}  {\|u\|}_p^p\!=\!0$.

From \cite{HCl}, we know that solutions of \eqref{eq1.1} belong to $C_{loc}^{1,\alpha}(\R^N)$ for some $\alpha\in(0,1)$ and decay exponentially at infinity. For weak solution $u \!\in\! X $ to \eqref{eq1.1}, we can follow the smooth function approximation procedure in \cite{MgV} to prove $P(u)\!=\!0$.

\end{proof}

Of course, the functional $I$ in \eqref{energyF} is  well defined in the entire $X$, indeed $2<p<\min\{2^*,q^*\}$ in the assumptions of Theorems \ref{th1.1}-\ref{th1.4}.

The proof of the $C^1$ regularity of $I$ in $X$ is almost standard, see for instance \cite{lByF}.
In turn, $I': X\to X'$ is given by
\begin{equation}\label{I'}
I'(u)\phi=\int_{\mathbb{R}^N} \nabla u \nabla \phi dx+\int_{\mathbb{R}^N}|\nabla u|^{q-2}
\nabla u \nabla\phi dx -\int_{\mathbb{R}^N}|u|^{p-2}u\phi dx
\end{equation}
for all $u, \phi\in X$.
Note that \eqref{eq1.1} has always the trivial solution $0$, but \eqref{eq1.2} prevents the case to occur.

Now, we report for completeness a part of Lemma 4.1.2 in \cite{cpy}, which will be useful in Section \ref{cri}.

\begin{lemma}\label{ineq}\textnormal{(\textnormal{\cite[Lemma 4.1.2]{cpy}})}
Suppose $\gamma>1$. Then, for any $a>0$, $b>0$ there exists $C_{\gamma}=C(\gamma)>0$ such that
\[
(a+b)^{\gamma}\leq
\left\{ {\begin{array}{*{20}{c}}
a^{\gamma}+b^{\gamma}+C_{\gamma} ( a^{ \frac{\gamma}{2} } b^{ \frac{\gamma}{2} } ),~~~~\gamma\in(1,2],    \\
a^{\gamma}+b^{\gamma}+C_{\gamma}( a^{\gamma-1} b+a b^{\gamma-1}  ),~~~~
\gamma >2.
\end{array}}
\right.\]
\end{lemma}

We end this section with the following useful result.

\begin{lemma}\label{lem3}\textnormal{(\textnormal{\cite[Lemma 2.7]{LIMa}})}
Assume $s>1$. Let $\Omega$ be an open set in $\mathbb{R}^N$, $\alpha, \beta$ positive numbers and $a(x,\xi)$ in
$C(\Omega\times \mathbb{R}^N,\mathbb{R}^N)$  such that
\renewcommand{\theenumii}{\roman{enumii}}
\begin{enumerate}
\item $\alpha|\xi|^{s}\le a(x,\xi)\xi$ for all $(x,\xi)\in\Omega\times\mathbb{R}^N$,
\item $|a(x,\xi)|\le\beta|\xi|^{s-1}$ for all $(x,\xi)\in\Omega\times\mathbb{R}^N$,
\item $(a(x,\xi)-a(x,\eta))(\xi,\eta)>0$ for all $(x,\xi)\in\Omega\times\mathbb{R}^N$ with $\xi\ne\eta$,
\item $a(x, \gamma \xi)=\gamma |\gamma|^{p-2} a(x, \xi)$ for all $(x,\xi)\in\Omega\times\mathbb{R}^N$ and $\gamma\in\mathbb R\setminus\{0\}$.
\end{enumerate}
Consider $(u_{n})_{n}, u\in W^{1,s}(\Omega)$, then $\nabla u_{n}\to \nabla u$ in $L^{s}(\Omega)$ if and only if
$$\lim_{n\to\infty}\int_{\Omega}\Big(a(x,\nabla u_{n}(x))-a(x,\nabla u(x))\Big)\Big(\nabla u_{n}(x)-\nabla u(x)\Big)dx=0.$$
\end{lemma}

\section{$L^2$-subcritical case}\label{sub}
In this section, we focus on the $L^2$-subcritical case and obtain a global minimizer for $$ m(c)\!:=\!\inf _{u \in S_c} I(u),$$
where $S_{c}\!:=\!\Big\{ u \!\in \! X \!:=\! H^{1}({\mathbb{R}^N}) \!\cap\! D^{1,q}({\mathbb{R}^N}): {\|u\|}_2\!=\!c  \Big\}$. We start with the following lemma.

\begin{lemma} \label{2.1.} Let $c\!>\!0$ and assume \eqref{h1.1a}, \eqref{h1.1b}. Then
$$-\infty< m(c)\!:=\!\inf _{u \in S_c} I(u)<0.$$
\end{lemma}

\begin{proof}
Observing that
\begin{equation*}
\left\{\begin{array}{l}
p\delta_{p}<q(1+\delta_{q})<2~~~~\mbox{if}~~~~\frac{2N}{N+2}
\!<\!q\!<\!2\!<\!p\!<\!q(1+\frac{2}{N}), \\
p\delta_{p}<2<q(1+\delta_{q})~~~~\mbox{if}~~~~2\!<\!q\!<\!N~~~~\mbox{and}
~~~~2\!<\!p\!<\!2(1\!+\!\frac{2}{N}),
\end{array}\right.
\end{equation*}
so we have
$$
   p\delta_{p}<\min\{ 2, q(1+\delta_{q}) \}.
$$
For any fixed $u \in S_{c}$, the Gagliardo-Nirenberg inequality \eqref{equ2.2} indicates that
\begin{align} \label{lower-bounded1.1}
I(u)=\frac{1}{2} \|\nabla u\|_{2}^2+\frac{1}{q} \|\nabla u\|_{q}^q-\frac{1}{p} {\|u\|}_p^p
\geq \frac{1}{2} \|\nabla u\|_{2}^2-\frac{\mathcal{C}^p_{N,p}}{p} c^{p(1-\delta_p)}\left\|\nabla u\right\|_{2}^{p\delta_p}.
\end{align} 
From \eqref{lower-bounded1.1}, we obtain $m(c)>-\infty$. Next, we show that $m(c)\!<\!0$. Let $u \in S_{c}$ be fixed, then we have $u_{t}(x)\!=\!t^{\frac{N}{2}}u(tx)\in S_{c}$ and
\begin{equation}\label{3.1bi}
I(u_t)=\frac{t^{2}}{2} \|\nabla u\|_{2}^2+\frac{t^{q(1+\delta_{q})}}{q} \|\nabla u\|_{q}^q-\frac{t^{ p\delta_{p}}}{p} {\|u\|}_p^p<0
\end{equation}
for $t\!>\!0$ sufficiently small. This implies that $m(c)\!\leq\!I(u_t)\!<\!0$.
\end{proof}

\begin{remark} \label{2.1.1-2} Under some extra condition on $c$, we can relax the restriction on $p$ and $q$ in Lemma \ref{2.1.}. In fact, if $N\!\geq\!2$, $1\!<\!q\!<\!2\!<\!p\!<\!2(1\!+\!\frac{2}{N})$ and    $c\!>\!c^{*}:=\inf\{c\in(0,+\infty):~~m(c)\!<\!0\}$, we still have
$$
-\infty< m(c)\!:=\!\inf _{u \in S_c} I(u)<0.
$$
The reason is that, for any $u\in S_1$ and $c>0$, we have $u(c^{-\frac{2}{N}}x) \in S_{c}$ and
\begin{align*} 
m(c)\leq I(u(c^{-\frac{2}{N}}x))=\frac{c^{ \frac{2(N-2)}{N} } }{2}\|\nabla u \|_{2}^2+\frac{c^{ \frac{2(N-q)}{N} } }{q}\|\nabla u \|_{q}^q-\frac{c^2}{p} {\|u \|}_p^p.
\end{align*}
Since $\frac{2(N-q)}{N}\!<\!2$ and $\frac{2(N-2)}{N}\!<\!2$, we see that
$m(c)<0$ when $c>0$ is sufficiently large. Then the set $\{c\in(0,+\infty):~~m(c)\!<\!0\}\not=\emptyset$. Define  $$c^{*}:=\inf\{c\in(0,+\infty):~~m(c)\!<\!0\},$$
we have $m(c)<0$ provided $c\!>\!c^{*}$. The lower bound $m(c)>-\infty$ can be proved in a standard way as that of Lemma \ref{2.1.}. 
\end{remark}

\begin{lemma}\label{LemA3.5}
Let $c\!>\!0$ and assume \eqref{h1.1a}, \eqref{h1.1b}. Then  \\
\indent (i) The map $c\mapsto m(c)$ is continuous;\\
\indent (ii) If $c_1\!\in\!(0,c)$ and $c_2\!=\!\sqrt{c^2-c^2_1}$, we have $ m(c)< m(c_1)+m(c_2)$.
\end{lemma}

\begin{proof}
(i)~~Let $c>0$ and $\{c_n\}\subset(0,+\infty)$ such that $c_n \to c$, it is sufficient to prove that $m(c_n) \to m(c)$. For every $n\in \mathbb{N}^+$, there exists $u_n \in S_{c_n}$ such that $m(c_n) \leq I(u_n) <m(c_n)+1/n$. We first show that $\{u_{n}\}$ is bounded in $X$. In fact, for $n$ sufficiently large, we deduce from Lemma \ref{2.1.} that $m(c_n)\leq0$. In a fashion similar to \eqref{lower-bounded1.1}, we get
\begin{align} \label{3.2boud2}
0\geq I(u_n)=\frac{1}{2} \|\nabla u_n\|_{2}^2+\frac{1}{q} \|\nabla u_n\|_{q}^q-\frac{1}{p} {\|u_n\|}_p^p \geq \frac{1}{2} \|\nabla u_n\|_{2}^2-\frac{\mathcal{C}^p_{N,p}}{p} c_n^{p(1-\delta_p)}\left\|\nabla u_n\right\|_{2}^{p\delta_p}.
\end{align}
Then inequality \eqref{3.2boud2} gives
\begin{equation}\label{3.2bbb}
\|\nabla u_n\|_{2}\leq \biggl[\frac{2\mathcal{C}^p_{N,p}}{p} \biggr]^{\frac{1}{2-p\delta_p}}c_n^{\frac{p(1-\delta_p)}{2-p\delta_p}}\leq \biggl[\frac{2\mathcal{C}^p_{N,p}}{p} \biggr]^{\frac{1}{2-p\delta_p}}c^{\frac{p(1-\delta_p)}{2-p\delta_p}}+o_n(1).
\end{equation}
With this upper bound, we derive from \eqref{3.2boud2} and the Gagliardo-Nirenberg inequality \eqref{equ2.2} that
\begin{equation}\label{3.2bi}
 \frac{1}{q} \|\nabla u_n\|_{q}^q\leq\frac{1}{p} {\|u_n\|}_p^p
\leq \frac{\mathcal{C}^p_{N,p}}{p} c_n^{p(1-\delta_p)}\left\|\nabla u_n\right\|_{2}^{p\delta_p}\leq 2^{\frac{p\delta_p}{2-p\delta_p}}\biggl[\frac{\mathcal{C}^p_{N,p}}{p} \biggr]^{\frac{2}{2-p\delta_p}}c^{\frac{2p(1-\delta_p)}{2-p\delta_p}}+o_n(1).
\end{equation}
Therefore, $\{u_{n}\}$ is bounded in $X$. Now considering $v_n := \frac{c}{c_n}u_n \in S_c$, we have
\begin{align*}
m(c) \leq I(v_n)&=I(u_n)+\frac{1}{2}\biggl(\frac{c^2}{c^2_n}-1\biggr) \|\nabla u_n\|_{2}^2+\frac{1}{q}\biggl(\frac{c^q}{c^q_n}-1\biggr)\|\nabla u_n\|_{q}^q-\frac{1}{p}\biggl(\frac{c^p}{c^p_n}-1\biggr) {\|u_n\|}_p^p\\
&=I(u_n)+o_n(1),
\end{align*}
where we used the boundedness of $\{u_n\}$ and the fact that $c_n\to c$. Passing to the limit as $n \to +\infty$,
we deduce that
 $$  m(c) \leq \lim_{n \to +\infty} \inf m(c_n).  $$
In a similar way, let $\{w_n\}$ be a minimizing sequence for $m(c)$, which is also bounded, and
let $z_n := \frac{c_n}{c}w_n \in S_{c_n}$. Then we have
$$
m(c_n) \leq I(z_n) = I(w_n) + o_n(1) \Longrightarrow \lim_{n \to +\infty}\sup m(c_n) \leq m(c).   $$
(ii)~~For any fixed $c_1\in(0,c)$, we first claim that
\begin{equation} \label{CoTq4}
 m(\theta c_1)< \theta^2 m(c_1),~~~~\forall \theta>1.
\end{equation}
Let $\{u_n\}\subset S_{c_1}$ be a minimizing sequence for $m(c_1)$, then $ u_n(\theta^{-\frac{2}{N}}x) \in S_{\theta c_1}$. By using $\frac{2(N-q)}{N}\!<\!2$ and $\frac{2(N-2)}{N}\!<\!2$, for $\theta>1$, we have
\begin{align*}
m(\theta c_1)-\theta^2I(u_n)&\leq I(u_n(\theta^{-\frac{2}{N}}x))-\theta^2I(u_n)\\
&=\frac{\theta^{ \frac{2(N-2)}{N} }-\theta^2}{2}\|\nabla u_n\|_{2}^2+\frac{\theta^{ \frac{2(N-q)}{N} }-\theta^2}{q}\|\nabla u_n\|_{q}^q\leq0.
\end{align*}
As a consequence $m(\theta c_1)\leq\theta^2 m(c_1)$, with equality if and only if $\|\nabla u_n\|_{2}^2\to 0$ and $\|\nabla u_n\|_{q}^q\to0$ as $n\to\infty$. In view of these facts, inequality \eqref{equ2.2} indicates that $\|u_n\|_{p}^p\to0$. It must be that $m(\theta c_1)<\theta^2 m(c_1)$, otherwise, we obtain a contradiction
$$0 > m(c_1) = \lim_{n\to\infty} I(u_n)=  \frac{1}{2} \lim_{n\to\infty}\|\nabla u_n\|_{2}^2+\frac{1}{q} \lim_{n\to\infty}\|\nabla u_n\|_{q}^q-\frac{1}{p} \lim_{n\to\infty}{\|u_n\|}_p^p= 0.$$
In the same manner, we can get
\begin{equation} \label{CoTq5}
 m(\theta c_2)< \theta^2 m(c_2),~~~~\forall \theta>1.
\end{equation}
Finally, apply \eqref{CoTq4} with $\theta=\frac{c}{c_1}>1$ and \eqref{CoTq5} with $\theta=\frac{c}{c_2}>1$ respectively, we get
\begin{align*}
 m(c)\!=\!\frac{c^2_1}{c^2}m\biggl(\frac{c}{c_1}c_1\biggr) \!+\!\frac{c_2^2}{c^2}  m\biggl(\frac{c}{c_2}c_2\biggr)
\!< \!m(c_1)\!+\!m(c_2).
\end{align*}
\end{proof}

Applying Lemma \ref{LemA3.5}, we can prove the compactness of the minimizing sequences for $m(c)$.

\begin{lemma} \label{LeMa3.4.1}
Assume that $c\!>\!0$ and \eqref{h1.1a}, \eqref{h1.1b}. Let $\{w_{n}\} \subset S_c$ be a minimizing sequence for $m(c)$,
then $m(c)$ possesses another minimizing sequence $\{u_{n}\} \subset S_c$ such that
$$\left\|u_{n}-w_{n}\right\|_{X} \to 0,~~~~~~~~\big(I|_{S_c}\big)'\left(u_{n}\right) \to 0\qquad \text{as}~~~~n\to+\infty.$$
Moreover, $\{u_{n}\}$ is relatively compact in $X$ up to translations and hence $m(c)$ is attained. 
\end{lemma}

\begin{proof}
Since $\{w_{n}\} \subset S_c$ is a minimizing sequence of $m(c)$, from the Ekeland's variational principle (cf. \cite[Theorem2.4]{Mwlm}), we get a new minimizing sequence $\{u_{n}\} \subset S_c$ for $m(c)$ such that $\left\|u_{n}-w_{n}\right\|_{X} \to 0$, which is also a Palais-Smale
sequence for $I|_{S_c}$. Hence, we have $\big(I|_{S_c}\big)'\left(u_{n}\right) \to 0$.

In the same way as the proof of Lemma \ref{LemA3.5} (i), we obtain that $\{u_{n}\}$ is bounded in $X$. If $\mathop {\lim }\limits_{n  \to \infty} \sup _{y \in \mathbb{R}^{N}} \int_{B_{R}(y)}\left|u_{n}(x)\right|^{2}dx\!=\!0$ for any $R\!>\!0$, we can prove that $\left\|u_{n}\right\|_{r} \to 0$ for $2<r<2^{*}$ (see \cite[Lemma I.1]{LilP}). In particular, we have ${\|u_n\|}_p^p\to0$, this together with Lemma \ref{2.1.} lead to the following contradiction
\begin{align*}
0&>m(c)=\mathop {\lim }\limits_{n  \to \infty}I(u_n)=\frac{1}{2}\mathop {\lim }\limits_{n  \to \infty}{\|\nabla u_n\|}_2^2+\frac{1}{q} \mathop {\lim }\limits_{n  \to \infty}\|\nabla u_n\|_{q}^q-\frac{1}{p}\mathop {\lim }\limits_{n  \to \infty}{\|u_n\|}_p^p \geq0.
\end{align*}
Then, there exist an $\varepsilon_0>0$ and a sequence $\{y_{n}\} \subset \mathbb{R}^{N}$ such that
$$
\int_{B_{R}(y_n)}\left|u_{n}(x)\right|^{2} d x \geq \varepsilon_0>0
$$
for some $R>0$. Hence we have $u_{n}(x+y_n)\rightharpoonup u_c\not\equiv0$ in $X$ for some $u_c\in X$. Let $v_n(x):=u_{n}(x+y_n)-u_c$, then we see that $v_{n}\rightharpoonup 0$ in $X$. Therefore, we get
\begin{equation}\begin{aligned} \label{NlaG1} &\left\|u_{n}\right\|_{2}^2=\left\|u_{n}(\cdot+y_n)\right\|_{2}^2=\left\|v_{n}\right\|_{2}^2
+\left\|u_c\right\|_{2}^2+o_n(1), \\
&\left\|\nabla u_{n}\right\|_{2}^2=\left\|\nabla u_{n}(\cdot+y_n)\right\|_{2}^2=\left\|\nabla v_{n}\right\|_{2}^2
+\left\|\nabla u_c\right\|_{2}^2+o_n(1).
\end{aligned}\end{equation}
Moreover, by the compactness of the embedding of $D^{1,q}(\mathbb R^N)$ in $L^s_{loc}(\mathbb R^N)$ for any $1\le s<q^*$ we have
$v_n\to 0$ in $L^{s}(\omega)$ with $\omega\Subset\mathbb{R}^N$, $1\le s< q^*$. So, by using an increasing sequence of compact sets whose union is $\mathbb R^N$ and a diagonal argument, we get $v_n(x)\to 0$ for a.e. $x\in\mathbb R^N$.
From the Br\'{e}zis-Lieb lemma \cite{a7}, we have
$$ \left\|u_{n}\right\|_{p}^p=\left\|u_{n}(\cdot+y_n)\right\|_{p}^p=\left\|v_{n}\right\|_{p}^p
+\left\|u_c\right\|_{p}^p+o_n(1).$$
By a standard truncation argument (cf. \cite{aLyM}, \cite{ByF2}, \cite{Wlzk}) in what follows we prove that $\nabla u_{n} \rightarrow \nabla u_c $ a.e. on ${\R}^N$, up to subsequences.
Choose $\psi\in C_0^\infty(\mathbb R^N)$ such that
$0\le\psi\le 1$ in $\mathbb R^N$, $\psi(x)=1$ for every $x\in B_1(0)$ and $\psi(x)=0$ for every $x\in \mathbb R^N\setminus B_2(0)$. Now, take $R>0$ and define $\psi_R(x)=\psi(x/R)$ for $x\in\mathbb R^N$.
Recalling the definition of $I'(u)\phi$ in \eqref{I'} with $u=u_n$ and $\phi=(u_n-u_c)\psi_R$, we get
\begin{equation}\label{qog}\begin{aligned}
\int_{\mathbb R^N} &\biggl[\nabla u_n-\nabla u_c+|\nabla u_n|^{q-2} \nabla u_n-|\nabla u_c|^{q-2} \nabla u_c\biggr](\nabla u_n-\nabla u_c)\psi_R dx=\\
&I'(u_n)((u_n-u_c)\psi_R)-\int_{\mathbb R^N}\nabla u_n u_n \nabla \psi_R  dx-\int_{\mathbb R^N} |  \nabla u_n|^{q-2} \nabla u_n u_n  \nabla \psi_R dx\\
&+ \int_{\mathbb R^N} |u_n|^p \psi_R dx+\int_{\mathbb R^N} \nabla u_n u_c \nabla \psi_R+\int_{\mathbb R^N} |\nabla u_n|^{q-2}\nabla u_n u_c \nabla \psi_R dx\\
&-\int_{\mathbb R^N} |u_n|^{p-2}u_n u_c \psi_R dx-\int_{\mathbb R^N} \nabla u_n \nabla u_c \psi_R dx-\int_{\mathbb R^N} |\nabla u_c|^{q-2}\nabla u_c \nabla u_n \psi_R dx\\
&+\int_{\mathbb R^N} | \nabla u_c|^{q}\psi_R dx+\int_{\mathbb R^N} | \nabla u_c|^{2}\psi_R dx.
\end{aligned}\end{equation}
Since $\{u_n\}\subset S_c$ and $\big(I|_{S_c}\big)'\left(u_{n}\right) \to 0$, we have $I'(u_n)((u_n-u_c)\psi_R)\to0$ as $n\to\infty$. Moreover, by the definition of $\psi_R$ and the boundedness of $u_n, u_c\in X$, then
$$\int_{\mathbb R^N}\nabla u_n 
u_n\nabla \psi_R  dx, \int_{\mathbb R^N} \nabla u_n u_c \nabla \psi_R\to 0, \qquad n\to\infty,$$
$$\int_{\mathbb R^N} |\nabla u_n|^{q-2}\nabla u_n u_n \nabla \psi_R dx, \int_{\mathbb R^N} |\nabla u_n|^{q-2}\nabla u_n u_c \nabla \psi_R dx\to 0, \qquad n\to\infty.$$
Following \cite{ByF2}, it holds
$$\int_{\mathbb R^N} |\nabla u_c|^{q-2}\nabla u_c \nabla u_n \psi_R dx\to \int_{\mathbb R^N} |\nabla u_c|^{q} \psi_R dx, \qquad    n\to\infty,$$
$$\int_{\mathbb R^N} \nabla u_n \nabla u_c \psi_R dx\to \int_{\mathbb R^N} |\nabla u_c|^2 \psi_R dx, \qquad n\to\infty,$$
$$\int_{\mathbb R^N} |u_n|^p \psi_R dx, \int_{\mathbb R^N} |u_n|^{p-2}u_n u_c \psi_R dx\to \int_{\mathbb R^N} |u_c|^{p}\psi_R dx, \qquad n\to\infty.$$
So that, \eqref{qog} becomes, for $n\to\infty$ the following
$$\lim_{n\to\infty}\int_{\mathbb R^N} \biggl[\nabla u_n-\nabla u_c+|\nabla u_n|^{q-2}\nabla u_n-|\nabla u_c|^{q-2}\nabla u_c\biggr](\nabla u_n-\nabla u_c)\psi_R dx=0.$$
By virtue of Lemma \ref{lem3} applied to $a(x,\xi)=|\xi|^{s-2}\xi$  with $s=2$ and $s=q$, following all steps in \cite{ByF2} with $B_{2R}$ instead of $A_\varepsilon$, we prove $\nabla u_{n} \rightarrow \nabla u_c $ a.e. on ${\R}^N$, up to subsequences.

Now, applying Br\'{e}zis-Lieb lemma \cite{a7} again, we obtain
\begin{align} \label{NlaG1.2}
\left\|\nabla u_{n}\right\|_{q}^q=\left\|\nabla u_{n}(\cdot+y_n)\right\|_{q}^q=\left\|\nabla v_{n}\right\|_{q}^q
+\left\|\nabla u_c\right\|_{q}^q+o_n(1).
\end{align}
We next claim that $u_{n}(x+y_n)\rightarrow u_c\not\equiv0~~~~\mbox{in}~~~~L^2(\R^N)$, or equivalently $v_{n}\!\rightarrow \!0$ in $L^2(\R^N)$. Denote $\left\|u_c\right\|_{2}\!=\!c_1$. If $c_1\!=\!c$, the proof is completed by \eqref{NlaG1}. If $c_1\!<\!c$, we learn from \eqref{NlaG1}-\eqref{NlaG1.2} that
$$m(c)=I\left(u_{n}\right)+o_n(1)=I\big(u_{n}(\cdot+y_n)\big)+o_n(1)
=I\left(v_{n} \right)+I\left(u_{c} \right)+o_n(1)\geq m(\left\| v_{n}\right\|_{2})+m(c_1).$$
By the continuity of $c\mapsto m(c)$ (see Lemma \ref{LemA3.5} (i)), we have
\begin{align} \label{NlaG4}
m(c)\geq m(c_2)+m(c_1),
\end{align}
where $c_2\!=\!\sqrt{c^2-c^2_1}\!>\!0$. However, \eqref{NlaG4} contradicts to  Lemma \ref{LemA3.5} (ii). Therefore, we have $\left\|u_c\right\|_{2}\!=\!c$ and hence $v_{n}\!\rightarrow \!0$ in $L^2(\R^N)$. It follows immediately, by \eqref{equ2.2}, that
$$\|v_n\|_{p} \leq \mathcal{C}_{N,p} \left\|\nabla v_n\right\|_{2}^{\delta_p} \left\|v_n\right\|_{2}^{(1-\delta_p)} \to0.$$
Finally, we get 
\begin{align*}
m(c)=I\left(v_{n} \right)+I\left(u_{c} \right)+o_n(1)\geq\frac{1}{2}{\|\nabla v_{n}\|}_2^2+\frac{1}{q} \|\nabla v_n\|_{q}^q+m(c)+o_n(1),
\end{align*}
which indicates ${\|\nabla v_{n}\|}_2\leq o_n(1)$ and ${\|\nabla v_{n}\|}_q\leq o_n(1)$. So we have $v_{n}\!\rightarrow \!0$ in $X$ and $u_{n}(\cdot+y_n)\rightarrow u_c\not\equiv0$ in $X$.
\end{proof}

\begin{proof}[\textbf {Proof of Theorem \ref{th1.1}.} ]
By using Lemma \ref{LeMa3.4.1}, we see that
$m(c)$ is attained by some $u_c \!\in\! S_c $. From \cite [Section 3.3 and Lemma 7.17]{ELMA}, we have $I\left({|u_c|}^*\right)\!\leq\!I\left(|u_c|\right)\!\leq\!I\left(u_c\right)$ since
$$\|\nabla{|u|}^*\|_2  \leq  \|\nabla{|u|} \|_2 ,~~~~\|\nabla{|u|}^*\|_q  \leq  \|\nabla{|u|} \|_q ,~~~~\|u\|_p=\|{|u|}^*\|_p,$$
where ${|u|}^*$ is the symmetric decreasing rearrangement of $|u|$. So we can assume that $u_c\!\in\! S_c$ is nonnegative and radially decreasing. Next, the Lagrange multipliers rule implies the existence of some $\lambda_c\in \mathbb{R}$ such that
\begin{align} \label{3.9eq4.2.1}
\int_{{\R^N}}\nabla u_c\nabla\varphi \!+\!\int_{\mathbb{R}^{N}} {| \nabla u_c |^{q - 2}}\nabla u_c \nabla{\varphi}\!-\! \int_{{\mathbb{R}^N}} \left|u_c\right|^{p\!-\!2} u_c {\varphi}\!-\!\lambda_c \int_{{\R^N}} u_c{\varphi}=0,~~~~\forall \varphi \in X.
\end{align}
That is, $(\lambda_c,u_c)$ satisfies \eqref{eq1.1}. From \eqref{3.9eq4.2.1} and Lemma   \ref{lem2.7}, we get
\begin{equation}\label{3.8bi}
\lambda_cc^2=\|\nabla u_c\|_{2}^2+\|\nabla u_c\|_{q}^q-{\|u_c\|}_p^p=\biggl(1-\frac{1}{\delta_{p}}\biggr)\|\nabla u_c\|_{2}^2+\biggl(1-\frac{ 1+\delta_{q} }{\delta_{p}}\biggr)\|\nabla u_c\|_{q}^q<0.\end{equation}
In a manner similar to get \eqref{3.2bbb}, we obtain 
$$\|\nabla u_c\|_{2}\leq \biggl[\frac{2\mathcal{C}^p_{N,p}}{p} \biggr]^{\frac{1}{2-p\delta_p}}c^{\frac{p(1-\delta_p)}{2-p\delta_p}}$$ 
Furthermore, as for \eqref{3.2bi}, we derive 
$$ \frac{1}{q} \|\nabla u_c \|_{q}^q\leq 2^{\frac{p\delta_p}{2-p\delta_p}}\biggl[\frac{\mathcal{C}^p_{N,p}}{p} \biggr]^{\frac{2}{2-p\delta_p}}c^{\frac{2p(1-\delta_p)}{2-p\delta_p}}. $$
As $u_c \in \mathcal{P}_{c}$, we can rewrite $I(u_c)$, by using Lemma \ref{lem2.7}, in another form
$$\begin{aligned}
0&>m(c)=I(u_c)=\frac{1}{2} \|\nabla u_c \|_{2}^2+\frac{1}{q} \|\nabla u_c \|_{q}^q-\frac{1}{p} {\| u_c \|}_p^p  \nonumber \\
&=\biggl(\frac{1}{2}-\frac{1}{p\delta_{p}}\biggr) \|\nabla u_c \|_{2}^2+\biggl(\frac{1}{q}-\frac{1+\delta_{q}}{p\delta_{p}}\biggr)\|\nabla u_c \|_{q}^q \nonumber \\
&\geq \biggl(\frac{1}{2}-\frac{1}{p\delta_{p}}\biggr) \biggl[\frac{2\mathcal{C}^p_{N,p}}{p} \biggr]^{\frac{2}{2-p\delta_p}}c^{\frac{2p(1-\delta_p)}{2-p\delta_p}}+
\biggl(1-\frac{q(1+\delta_{q})}{p\delta_{p}}\biggr)  2^{\frac{p\delta_p}{2-p\delta_p}}\biggl[\frac{\mathcal{C}^p_{N,p}}{p} \biggr]^{\frac{2}{2-p\delta_p}}c^{\frac{2p(1-\delta_p)}{2-p\delta_p}}.
\end{aligned}$$
Therefore, we have $m(c)\gtrsim-c^{\frac{2p(1-\delta_p)}{2-p\delta_p}}$. This fact also indicates that $m(c) \to 0^{-}$ as $c \to 0^{+}$. Similarly, from \eqref{3.8bi}, we have $\lambda_c \gtrsim -c^{\frac{2(p-2)}{2-p\delta_p}}$. This fact leads to $\lambda_c \to 0^{-}$ as $c \to 0^{+}$.
\end{proof}

\section{$L^2$ and $L^q$ critical case}\label{cri}
In this section, we consider the $L^2$ and $L^q$ critical cases, when we will obtain several nonexistence results. First, we state some preliminary results.

\begin{lemma} \label{2.1.-0} Let $N\!\geq\!2$, $p\!=\!\bar{p}\!:=\!2(1\!+\!\frac{2}{N})$, $m(c)\!=\!\inf _{u \in S_c} I(u)$ and $c_{*}:=\|W_{\bar{p}}\|_{2}$, with $W_{p}$ defined in Lemma \ref{lem2.1}. Then, we have  \\
(1) If $1\!<\!q\!<\!N$ and $0<c\leq c_{*}$, then $m(c)=0$ and there is no minimizer of $m(c)$;  \\
(2) If $1\!<\!q\!<\!2$ and $c>c_{*}$, then $m(c)=-\infty$ and there is no minimizer of $m(c)$.
\end{lemma}

\begin{proof}
(1) We first observe that $p\delta_{p}=2$ if $p\!=\!2(1\!+\!\frac{2}{N})$. For any fixed $u \in S_{c}$, the Gagliardo-Nirenberg inequality \eqref{equ2.2} indicates that
\begin{align} \label{4-lower-bounded1.1}
I(u)=\frac{1}{2} \|\nabla u\|_{2}^2+\frac{1}{q} \|\nabla u\|_{q}^q-\frac{1}{\bar{p}} {\|u\|}_{ \bar{p} }^{ \bar{p} }
\geq \biggl(\frac{1}{2}-\frac{\mathcal{C}^{ \bar{p} }_{N,\bar{p}}}{\bar{p}} c^{\bar{p}-2}\biggr) \left\|\nabla u\right\|_{2}^{2}+\frac{1}{q} \|\nabla u\|_{q}^q\geq0
\end{align} 
when $c \in (0, \|W_{\bar{p}}\|_{2}]$. So we have $m(c)\geq0$. In addition, we have $u_{t}(x)\!=\!t^{\frac{N}{2}}u(tx)\in S_{c}$ and hence
\begin{equation}\label{upperbound}
m(c)\!\leq\!I(u_t)=\frac{t^{2}}{2} \|\nabla u\|_{2}^2+\frac{ t^{q(1+\delta_{q})} }{q} \|\nabla u\|_{q}^q-\frac{t^{ 2}}{ \bar{p} } {\|u\|}_{ \bar{p} }^{ \bar{p} } \to 0
\end{equation}
as $t\!\to\!0$. This implies that $m(c)\!=\!0$. From \eqref{4-lower-bounded1.1}, we also deduce that $m(c)$ can not be attained.

(2) First, note that arguing as in the proof of Lemma \ref{lem2.7}, we get
\begin{equation} \label{TeStfunC3.3}
 \| \nabla W_p \|_{2}^{2}=\| W_p \|_{2}^{2}=\frac{2}{p} \| W_p \|_{p}^{p}.
\end{equation}
In addition, it follows from \cite[Proposition 4.1]{Gdas} that
\begin{equation} \label{TeStfunC3.4}
W_p(x),|\nabla W_p(x)|=O\left(|x|^{-\frac{N-1}{2}} e^{-|x|}\right)  \quad \text { as }|x| \rightarrow \infty.
\end{equation}
Let $\varphi(x) \in C_{0}^{\infty}\left(\mathbb{R}^{N}\right)$ be a nonnegative radial function such that $\varphi(x)=1$ if $|x| \leq 1$ and $\varphi(x)=0$ if $|x| \geq 2$. For $\tau>0$ and $R>0$, we denote
$$\phi_{1}(x)=A_{\tau, R} \frac{ (\tau c)^{ \frac{N}{2} }  }{\|W_p\|_{2}} \varphi \Big( \frac{x}{R} \Big) W_p ( \tau x ) ,$$
where $A_{\tau, R}>0$ is chosen such that $\int_{\mathbb{R}^{N}} \phi_{1}^{2} d x=c^2$, which imply
$$\begin{aligned}
c^2=\int_{\mathbb{R}^{N}} \phi_{1}^{2} d x=A^2_{\tau, R} \frac{ (\tau c)^{ N }  }{\|W_p\|^2_{2}} \int_{\mathbb{R}^{N}} \varphi^2 \Big( \frac{x}{R} \Big) W^2_p ( \tau x )dx   =A^2_{\tau, R} \frac{  c^{ N }}{\|W_p\|^2_{2}} \int_{\mathbb{R}^{N}} \varphi^2 \Big( \frac{y}{\tau R} \Big) W^2_p ( y )dy.
 \end{aligned}$$
By scaling and \eqref{TeStfunC3.4}, we have
\begin{align}   \label{aSTimateTEfunC}
\frac{1}{A_{\tau, R}^{2}}=&\frac{ c^{N-2} }{\| W_p \|_{2}^{2}} \int_{\mathbb{R}^{N}} \varphi^{2}\left(\frac{y}{\tau R}\right)  |W_p(y)|^{2} d y =\frac{ c^{N-2} }{\| W_p \|_{2}^{2}} \int_{\mathbb{R}^{N}} \Big( \varphi^{2}\left(\frac{y}{\tau R}\right)-1+1 \Big)  |W_p(y)|^{2} d y  \nonumber \\
=&c^{N-2}+\frac{ c^{N-2} }{\| W_p \|_{2}^{2}} \int_{\mathbb{R}^{N}} \Big( \varphi^{2}\left(\frac{y}{\tau R}\right)-1\Big)  |W_p(y)|^{2} d y  \nonumber \\
=&c^{N-2}+\frac{ c^{N-2} }{\| W_p \|_{2}^{2}} \Big\{ \int_{ B_{\tau R}(0) } \Big( \varphi^{2}\left(\frac{y}{\tau R}\right)-1\Big)  |W_p(y)|^{2}dy   \nonumber \\
&+\int_{ \mathbb{R}^{N} / B_{\tau R}(0) } \Big( \varphi^{2}\left(\frac{y}{\tau R}\right)-1\Big)  |W_p(y)|^{2}dy  \Big\}  \nonumber \\
=&c^{N-2}+\frac{ c^{N-2} }{\| W_p \|_{2}^{2}} \Big\{  \int_{ \mathbb{R}^{N} / B_{\tau R}(0) } \Big( \varphi^{2}\left(\frac{y}{\tau R}\right)-1\Big)  |W_p(y)|^{2}dy  \Big\}  \nonumber \\
=&c^{N-2}+\frac{ c^{N-2} }{\| W_p \|_{2}^{2}} \Big\{  \int_{ \mathbb{R}^{N} / B_{\tau R}(0) } \tilde{C}\Big( \varphi^{2}\left(\frac{y}{\tau R}\right)-1\Big)  |y|^{-(N-1)} e^{-2|y|}dy  \Big\}  \nonumber  \\
=&c^{N-2}+\frac{ c^{N-2} }{\| W_p \|_{2}^{2}} \Big\{  \int_{ \tau R }^{+\infty} \tilde{C}\Big( \varphi^{2}\left(\frac{r}{\tau R}\right)-1\Big)  r^{-(N-1)} e^{-2r} \omega_N r^{(N-1)} dr  \Big\}  \nonumber  \\
=&c^{N-2}+\frac{ c^{N-2} }{\| W_p \|_{2}^{2}} \Big\{  \int_{ \tau R }^{+\infty} \tilde{C}\omega_N \Big( \varphi^{2}\left(\frac{r}{\tau R}\right)-1\Big)  e^{-2r}  dr  \Big\}  \nonumber  \\
=&c^{N-2}+O\left((\tau R)^{-\infty}\right)
 \text { as } \tau R \rightarrow \infty,
\end{align}
where $\tilde{C}>0$ is a constant, $\omega_N=\int_{ \partial B_{1}(0) }dS$, $f(t)=O\left(t^{-\infty}\right)$ means $\lim _{t \rightarrow \infty}|f(t)| t^{s}=0$ for all $s>0$. Note that
$$\begin{aligned}    
\nabla \phi_{1}&=A_{\tau, R} \frac{ (\tau c)^{ \frac{N}{2} }  }{\|W_p\|_{2}} \biggl[R^{-1}\nabla\varphi \Big( \frac{x}{R} \Big) W_p ( \tau x )+\tau\varphi \Big( \frac{x}{R} \Big) \nabla W_p ( \tau x )\biggr]\\
\end{aligned}.$$
Taking, for instance, $R=\tau^{-\frac{1}{2}}$, then \eqref{TeStfunC3.3}-\eqref{aSTimateTEfunC} indicate that
\begin{align} \label{4.6aSTimateTEfunC3.7}
\int_{\mathbb{R}^{N}}\left|\nabla \phi_{1}\right|^{2} d x
 =& \frac{A_{ \tau, R }^{2} \tau^{N} c^N}{\|W_p\|_{2}^{2}}\Big \{ \frac{1}{\tau^{N} R^2} \int_{\mathbb{R}^{N}} \Bigl|\nabla \varphi\Bigl(\frac{x}{\tau R}\Bigr)\Bigr|^{2} | W_p(x)|^{2} d x   \nonumber \\
 &+  \tau^{2-N} \int_{\mathbb{R}^{N}} \Bigl| \varphi\Bigl(\frac{x}{\tau R}\Bigr)\Bigr|^{2} |\nabla W_p(x)|^{2} d x +O\left((\tau R)^{-\infty}\right) \Big \} \nonumber \\
=& \tau^{2} \frac{ c^{2} }{\|W_p\|_{2}^{2}}  \int_{\mathbb{R}^{N}} |\nabla W_p(x)|^{2} d x  +O\left((\tau R)^{-\infty}\right),
\end{align}


By using Lemma \ref{ineq} with $\gamma=q\in(1,2)$, we have the following estimates
\begin{align}
 \int_{\mathbb{R}^{N}}&\left|\nabla \phi_{1}\right|^{q} d x
 \leq  \frac{  A_{ \tau, R }^{q} \tau^{\frac{Nq}{2}} c^{\frac{Nq}{2}} }{\|W_p\|_{2}^{q}}   \int_{\mathbb{R}^{N}}  \Big \{ |  \frac{1}{R}\nabla\varphi \Big( \frac{x}{R} \Big) W_p ( \tau x ) |^{q}+| \tau\varphi \Big( \frac{x}{R} \Big) \nabla W_p ( \tau x )|^{q}   \nonumber \\
 &+C_{q} \Big( | \frac{1}{R}\nabla\varphi \Big( \frac{x}{R} \Big) W_p ( \tau x ) |^{ \frac{q}{2} } | \tau\varphi \Big( \frac{x}{R} \Big) \nabla W_p ( \tau x ) |^{ \frac{q}{2} } \Big)  \Big \} d x  \nonumber \\
 = & \frac{  A_{ \tau, R }^{q} \tau^{\frac{Nq}{2}} c^{\frac{Nq}{2}} }{\|W_p\|_{2}^{q}}    \Big \{   \frac{1}{\tau^{N} R^q} \int_{\mathbb{R}^{N}} |\nabla \varphi(\frac{x}{\tau R})|^{q} | W_p(x)|^{q} d x + \tau^{q-N} \int_{\mathbb{R}^{N}} | \varphi(\frac{x}{\tau R})|^{q} |\nabla W_p(x)|^{q} dx   \nonumber \\
 &+C_{q} \frac{ \tau^{ \frac{q}{2}-N }   } { R^{ \frac{q}{2} }  }  \int_{\mathbb{R}^{N}} |  \nabla\varphi \Big( \frac{x}{\tau R} \Big) W_p ( x ) |^{ \frac{q}{2} } | \varphi \Big( \frac{x}{\tau R} \Big) \nabla W_p ( x ) |^{ \frac{q}{2} } d x   \Big \}  \nonumber \\
 = & \frac{  A_{ \tau, R }^{q} \tau^{\frac{Nq}{2}} c^{\frac{Nq}{2}} }{\|W_p\|_{2}^{q}}    \Big \{   \frac{1}{\tau^{N} R^q} \int_{ \mathbb{R}^{N} / B_{\tau R}(0) } |\nabla \varphi(\frac{x}{\tau R})|^{q} | W_p(x)|^{q} d x    \nonumber \\
 &+ \tau^{q-N} \|\nabla W_p \|_{q}^{q}+ \tau^{q-N} \int_{ \mathbb{R}^{N} / B_{\tau R}(0) } \Big( | \varphi(\frac{x}{\tau R})|^{q}-1\Big) |\nabla W_p(x)|^{q} dx   \nonumber \\
 &+C_{q} \frac{ \tau^{ \frac{q}{2}-N }   } { R^{ \frac{q}{2} }  }  \int_{ \mathbb{R}^{N} / B_{\tau R}(0) } |  \nabla\varphi \Big( \frac{x}{\tau R} \Big) W_p ( x ) |^{ \frac{q}{2} } | \varphi \Big( \frac{x}{\tau R} \Big) \nabla W_p ( x ) |^{ \frac{q}{2} } d x   \Big \}  \nonumber \\
 = & \frac{  A_{ \tau, R }^{q} \tau^{\frac{Nq}{2}} c^{\frac{Nq}{2}} }{\|W_p\|_{2}^{q}}    \Big \{   \frac{1}{\tau^{N} R^q} \int_{ \mathbb{R}^{N} / B_{\tau R}(0) } C |\nabla \varphi(\frac{x}{\tau R})|^{q} |x|^{\frac{q(1-N)}{2}} e^{-q|x|} d x    \nonumber \\
 &+ \tau^{q-N} \|\nabla W_p \|_{q}^{q}+ \tau^{q-N} \int_{ \mathbb{R}^{N} / B_{\tau R}(0) } C \Big( | \varphi(\frac{x}{\tau R})|^{q}-1\Big) |x|^{\frac{q(1-N)}{2}} e^{-q|x|} dx   \nonumber \\
 &+C_{q} \frac{ \tau^{ \frac{q}{2}-N }   } { R^{ \frac{q}{2} }  }  \int_{ \mathbb{R}^{N} / B_{\tau R}(0) }C  |  \varphi \Big( \frac{x}{\tau R} \Big) \nabla\varphi \Big( \frac{x}{\tau R} \Big)   |^{ \frac{q}{2} }   |x|^{\frac{q(1-N)}{2}} e^{-q|x|} d x   \Big \}  \nonumber \\
 = & \frac{  A_{ \tau, R }^{q} \tau^{\frac{Nq}{2}} c^{\frac{Nq}{2}} }{\|W_p\|_{2}^{q}}    \Big \{   \frac{1}{\tau^{N} R^q}  \int_{\tau R}^{+\infty}  \tilde{C} |\nabla \varphi(\frac{r}{\tau R})|^{q} r^{\frac{(q-2)(1-N)}{2}} e^{-qr} d r    \nonumber \\
 &+ \tau^{q-N} \|\nabla W_p \|_{q}^{q}+ \tau^{q-N} \int_{\tau R}^{+\infty}  \tilde{C} \Big( | \varphi(\frac{r}{\tau R})|^{q}-1\Big) r^{\frac{(q-2)(1-N)}{2}} e^{-qr} dr   \nonumber \\
 &+C_{q} \frac{ \tau^{ \frac{q}{2}-N }   } { R^{ \frac{q}{2} }  }  \int_{\tau R}^{+\infty} \tilde{C}  |  \varphi \Big( \frac{r}{\tau R} \Big) \nabla\varphi \Big( \frac{r}{\tau R} \Big)   |^{ \frac{q}{2} }   r^{\frac{(q-2)(1-N)}{2}} e^{-qr} d r   \Big \}  \nonumber \\
 \leq & \frac{  A_{ \tau, R }^{q} \tau^{\frac{Nq}{2}} c^{\frac{Nq}{2}} }{\|W_p\|_{2}^{q}}    \Big \{   \tau^{q-N} \|\nabla W_p \|_{q}^{q}+ 
 O\left((\tau R)^{-\infty}\right)
 \Big \}  \nonumber \\
=& \tau^{ q(1+\delta_{q}) } \frac{ c^q}{\|W_p\|_{2}^{q}}  \int_{\mathbb{R}^{N}} |\nabla W_p(x)|^{q} d x  +O\left((\tau R)^{-\infty}\right).\nonumber 
\end{align}

So we conclude that 
\begin{align} \label{4.7aSTimateTEfunC3.8}
\int_{\mathbb{R}^{N}}\left|\nabla \phi_{1}\right|^{q} d x
 \leq  \tau^{ q(1+\delta_{q}) } \frac{  c^q}{\|W_p\|_{2}^{q}}  \int_{\mathbb{R}^{N}} |\nabla W_p(x)|^{q} d x  +O\left((\tau R)^{-\infty}\right).
\end{align}
Making slight modifications, we can see that it remains true when $q\in(2,N)$.
Moreover,
\begin{align} \label{4.8aSTimateTEfunC3.9}
\int_{\mathbb{R}^{N}}\left|\phi_{1}\right|^{p} d x
 = & \frac{  A_{ \tau, R }^{p} \tau^{ \frac{Np}{2}-N } c^{ \frac{Np}{2} } }{\|W_p\|_{2}^{p}}\int_{\mathbb{R}^{N}} \Bigl| \varphi\Bigl(\frac{x}{\tau R}\Bigr)\Bigr|^{p} | W_p(x)|^{p} d x  \nonumber \\
=& \tau^{ p\delta_{p} } \frac{c^p}{\|W_p\|_{2}^{p}}\int_{\mathbb{R}^{N}}  | W_p(x)|^{p} d x  +O\left((\tau R)^{-\infty}\right).
\end{align}
Letting $\tau \to +\infty$, it follows from \eqref{4.6aSTimateTEfunC3.7}-\eqref{4.8aSTimateTEfunC3.9} that
\begin{align} \label{lowerbound-4.9}
m(c) \leq & I(\phi_{1}) \leq \frac{\tau^{2}}{2} \frac{ c^2}{\|W_p\|_{2}^{2}}  \int_{\mathbb{R}^{N}} |\nabla W_p(x)|^{2} d x+\frac{ \tau^{q(1+\delta_{q})} }{q} \frac{ c^q}{\|W_p\|_{2}^{q}}  \int_{\mathbb{R}^{N}} |\nabla W_p(x)|^{q} d x  \nonumber \\
&-\frac{\tau^{ p\delta_{p} }}{p}  \frac{c^p}{\|W_p\|_{2}^{p}}\int_{\mathbb{R}^{N}}  | W_p(x)|^{p} d x  +O\left((\tau R)^{-\infty}\right), \nonumber \\
&= \tau^{2}c^2  \biggl(  \frac{ \|\nabla W_p\|_2^{2}}{2\|W_p\|_{2}^{2}}  -   \frac{ c^{p-2} \| W_p\|_p^{p}  }{p\|W_p\|_{2}^{p}}  \biggr)  +\frac{\tau^{ q(1+\delta_{q}) }}{q} \frac{  c^q}{\|W_p\|_{2}^{q}}  \|W_p\|_{q}^{q} +O\left((\tau R)^{-\infty}\right).
\end{align}
Since $1\!<\!q\!<\!2$ and $c>\|W_{\bar{p}}\|_{2}$, we have
$$ \frac{ \|\nabla W_{\bar{p}}\|_2^{2}}{2\|W_{\bar{p}}\|_{2}^{2}}  -   \frac{ c^{{\bar{p}}-2} \| W_{\bar{p}}\|_{\bar{p}}^{{\bar{p}}}  }{{\bar{p}}\|W_{\bar{p}}\|_{2}^{{\bar{p}}}}<0\text{ and } q(1+\delta_{q})<2 .$$
Letting $\tau \to +\infty$ in \eqref{lowerbound-4.9}, we have $m(c) =-\infty$.
\end{proof}

\begin{lemma} \label{2.1.-0bis} Let $N\!\geq\!3$, $p\!=\!\hat{p}:=\!q(1+\frac{2}{N})$, $m(c)\!:=\!\inf _{u \in S_c} I(u)$, $W_{p}$ be defined in Lemma \ref{lem2.1} and $\mathcal{K}_{N,p}$ in Lemma \ref{lem2.2}. Let
$$c_{**}:=\biggl(\frac{N+2}{N\mathcal{K}_{N,p}^{q(N+2)/N}}\biggr)^{N/2q}, \qquad \hat{c}_{**}:=\biggl[\frac{ 2\|\nabla W_p\|_q^q}{q\|W_p\|_2^{2(N-q)/N}}\biggr]^{N/2q}.$$
Then, we have  \\
(1) If $1\!<\!q\!<\!N$ and $0<c\leq c_{**}$, then $m(c)=0$ and there is no minimizer of $m(c)$;  \\
(2) If $2\!<\!q\!<\!N$ and $c>\hat{c}_{**}$, then $m(c)=-\infty$ and there is no minimizer of $m(c)$.
\end{lemma}

\begin{proof}
Following Lemma \ref{2.1.-0}, by using \eqref{pq-equ2.2} in place of \eqref{equ2.2} and observing that $p\nu_{p,q}=q$ if $p\!=\!q(1+\frac{2}{N})$, we have     
\begin{align*}
I(u)&=\frac{1}{2} \|\nabla u\|_{2}^2+\frac{1}{q} \|\nabla u\|_{q}^q-\frac{1}{ {p}} {\|u\|}_{  {p} }^{  {p} }  \\
&\geq \frac{1}{2} \|\nabla u\|_{2}^2+\frac{1}{q} \|\nabla u\|_{q}^q-\frac{\mathcal{K}^p_{N,p}}{ {p}} \left\|\nabla u\right\|_{q}^{p\nu_{p,q} } c^{p(1-\nu_{p,q} )}  \\
&=\frac{1}{2} \|\nabla u\|_{2}^2+ \Big(\frac{1}{q}-\frac{\mathcal{K}^p_{N,p}}{ {p}}  c^{p-q}  \Big) \left\|\nabla u\right\|_{q}^{q}  \geq 0
\end{align*}
when $\frac{\mathcal{K}^p_{N,p}}{ {p}}  c^{p-q}\le\frac{1}{q}$, that is if $c\le c_{**}$. So that, taking into account \eqref{upperbound}, if $0<c\le c_{**}$, then $m(c)=0$ and there is no minimizer of $m(c)$. Thus, (1) is proved.

Concerning (2), we can repeat the argument of Lemma \ref{2.1.-0} 
in the same way as in \eqref{4.6aSTimateTEfunC3.7}, from \eqref{4.7aSTimateTEfunC3.8} which holds also if when $q\in(2,N)$ by using Lemma \ref{ineq}, since $p\delta_p=q(1+\delta_q)$ then \eqref{lowerbound-4.9} becomes
$$m(c)\leq \frac{\tau^{2}}{2} \frac{ c^2}{\|W_p\|_{2}^{2}}  \int_{\mathbb{R}^{N}} |\nabla W_p(x)|^{2} d x+\tau^{ p\delta_{p} }c^q\biggl[\frac{ \|\nabla W_p\|_q^q}{q\|W_p\|_{2}^{q}}-\frac{c^{p-q}\|W_p\|_p^p}{p\|W_p\|_{2}^{p}}\biggr]+O\left((\tau R)^{-\infty}\right).$$  
Note that $2\!<\!q\!<\!N$, so we have $ p\delta_{p}=q(1+\delta_{q})>2$. In view of this fact and $c>\hat{c}_{**}$, we let $\tau \to +\infty$ and obtain $m(c) =-\infty$.  



\end{proof}

\begin{remark}\label{rem4.2}
By using \eqref{pq-equ2.2}, we can easily see that $c_{**}<\hat{c}_{**}$.
Thus, by the arguments above, we get nonexistence for $c<c_{**}$ and $c>\hat{c}_{**}$. The gap in $[c_{**},\hat{c}_{**}]$ is still an open problem and it is caused by the test function. To obtain the optimal results, we shall use the externals of \eqref{pq-equ2.2}, denoted by $W_{p,q}$. All we need is the exponential decay of $W_{p,q}$. In other words, if 
$W_{p,q}(x),|\nabla W_{p,q}(x)|=O\left(|x|^{-\frac{N-1}{2}} e^{-|x|}\right)$ as $|x| \rightarrow \infty$. Indeed, if it holds, thus, letting
$$\psi_{1}(x)=A_{\tau, R} \frac{ (\tau c)^{ \frac{N}{2} }  }{\|W_{p,q}\|_{2}} \varphi \Big( \frac{x}{R} \Big) W_{p,q} ( \tau x ) ,$$ 
as a test function, then $\hat{c}_{**}$ should be replaced by $c_{**}$ and the result will become sharp.

\end{remark}

\begin{proof}[\textbf {Proofs of Theorems \ref{th1.2} and \ref{th1.2bis}} ]
It is sufficient to follow Lemma \ref{2.1.-0} and Lemma \ref{2.1.-0bis}, respectively.
\end{proof}

\section{$L^2$-supercritical case}\label{supc}
In this section, we consider the $L^2$-supercritical case. Subsection \ref{subsec5.1} is devoted to proving the existence of a ground state, that is Theorem 1.3, whose statement is given in the Introduction. We also prove Theorem 1.4 which implies the existence of infinitely many radial solutions in subsection \ref{subsec5.2}.

\subsection{Existence of a ground state} \label{subsec5.1}

We first observe that
\begin{equation}\label{cond3}\left\{\begin{array}{l}
q(1+\delta_{q})<2<p\delta_{p}~~\mbox{if}~~1\!<\!q\!<\!2~~\mbox{and}~~2(1\!+\!\frac{2}{N})\!<\!p\!<\!q^*, \\
2<q(1+\delta_{q})<p\delta_{p}~~\mbox{if}~~2\!<\!q\!<\!N~~\mbox{and}~~q(1+\frac{2}{N})\!<\!p\!<\!2^*.
\end{array}\right.
\end{equation}
Let $u \in S_{c}$ be fixed, then we have $u_{t}(x)\!=\!t^{\frac{N}{2}}u(tx)\in S_{c}$ and
\begin{equation*}
m(c)\!=\!\inf _{u \in S_c}  I(u)\!\leq\!I(u_t)=\frac{t^{2}}{2} \|\nabla u\|_{2}^2+\frac{t^{q(1+\delta_{q})}}{q} \|\nabla u\|_{q}^q-\frac{t^{ p\delta_{p}}}{p} {\|u\|}_p^p \to -\infty~~~~\mbox{as}~~~~t \to +\infty.
\end{equation*}
In this case, $I(u)$ is unbounded from below on $S_c$, so the global minimization method is not valid in searching critical points of $I|_{S_c}$. For this reason, we recall the Pohozaev set:
$$\mathcal{P}_{c}=\left\{u \in S_{c} : P(u)\!:=\! \|\nabla u\|_{2}^2+(1+\delta_{q})  \|\nabla u\|_{q}^q-\delta_{p}  {\|u\|}_p^p\!=\!0 \right\}.$$
Actually, it is induced by $\frac{dI(u_{t})}{dt}|_{t=1}\!:=\!P(u)=0$ for $u_{t}(x)\!=\!t^{\frac{N}{2}}u(tx)$. We consider a modified minimization problem
$$ \sigma(c)\!:=\!\inf _{ u \in \mathcal{P}_{c} } I(u).$$
From the following lemma, we see that $I(u)$ is bounded from below on $\mathcal{P}_{c}$.

\begin{lemma} \label{pqlalpl3.1}
Let $c>0$ and assume \eqref{h1.3a}, \eqref{h1.3b}. Then, $I$ is coercive on $\mathcal{P}_{c}$ and
$$ \sigma(c)\!:=\!\inf _{u \in \mathcal{P}_{c}} I(u)>0.$$
\end{lemma}

\begin{proof}
For any $u \in \mathcal{P}_{c}$, we have $\|\nabla u\|_{2}^2+(1+\delta_{q})  \|\nabla u\|_{q}^q=\delta_{p}  {\|u\|}_p^p$. To begin with, we deduce from inequalities \eqref{equ2.2}- \eqref{pq-equ2.2} that
$$\begin{aligned}
\|\nabla u\|_{2}^2\leq\|\nabla u\|_{2}^2+(1+\delta_{q})  \|\nabla u\|_{q}^q=\delta_{p}  {\|u\|}_p^p\leq\delta_{p} \mathcal{C}^p_{N,p} \left\|\nabla u\right\|_{2}^{p\delta_p} c^{p(1-\delta_p)}, \nonumber \\
(1+\delta_{q})  \|\nabla u\|_{q}^q\leq\|\nabla u\|_{2}^2+(1+\delta_{q})  \|\nabla u\|_{q}^q=\delta_{p}  {\|u\|}_p^p\leq\delta_{p} \mathcal{K}^p_{N,p} \left\|\nabla u\right\|_{q}^{p\nu_{p,q}} c^{p(1-\nu_{p,q})},
\end{aligned}$$
which lead to
\begin{align} \label{pq4.3.1}
\left\|\nabla u\right\|_{2}\geq\biggl[\frac{1}{\delta_{p} \mathcal{C}^p_{N,p}}\biggr]^{ \frac{1}{p\delta_{p}-2} }c^{-\frac{p(1-\delta_{p})}{p\delta_{p}-2}}~~\mbox{and}~~
\left\|\nabla u\right\|_{q}\geq\biggl[\frac{1+\delta_{q}}{\delta_{p} \mathcal{K}^p_{N,p}}\biggr]^{ \frac{1}{p\nu_{p,q}-q} }c^{-\frac{p(1-\nu_{p,q})}{p\nu_{p,q}-q}}.
\end{align}
Furthermore, for each $u \in \mathcal{P}_{c}$, we can rewrite $I(u)$ in another form
\begin{align} \label{pq4.4}
I(u)=\frac{1}{2} \|\nabla u\|_{2}^2+\frac{1}{q} \|\nabla u\|_{q}^q-\frac{1}{p} {\|u\|}_p^p=\biggl(\frac{1}{2}-\frac{1}{p\delta_{p}}\biggr) \|\nabla u\|_{2}^2+\biggl(\frac{1}{q}-\frac{1+\delta_{q}}{p\delta_{p}}\biggr)\|\nabla u\|_{q}^q.
\end{align}
Consequently, for every sequence $\{u_k\}\subset \mathcal{P}_{c}$ such that ${\|u_k\|}_{X}\to+\infty$, we deduce from $p\delta_{p}>2$ and $p\delta_{p}>q(1+\delta_{q})$ that $I(u_k)\to+\infty$, hence we have $I$ is coercive on $\mathcal{P}_{c}$. By using \eqref{pq4.3.1}-\eqref{pq4.4}, we also obtain  $\sigma(c)\!:=\!\inf _{ u \in \mathcal{P}_{c} } I(u)>0$.
\end{proof}

\begin{lemma} \label{3.2pqlalplace} Let $c>0$ and assume \eqref{h1.3a}, \eqref{h1.3b}. Then, for any $u\in S_c$ and $u_{t}(x)\!=\!t^{\frac{N}{2}}u(tx)$, there exists a unique $t_0 > 0$ such that $I(u_{t_0} ) = \max_{\{t>0\}}I(u_t)$ and $u_{t_0}\in \mathcal{P}_{c}$. In particular, we have
\begin{itemize}
\item[(i)] $t_0< 1\Leftrightarrow P(u)<0$,
\item[(ii)] $t_0=1 \Leftrightarrow P(u)=0$.
\end{itemize}
\end{lemma}

\begin{proof}
Following \eqref{3.1bi}, for any $u\in S_c$, let
$$h(t) :=I(u_t)=\frac{t^{2}}{2} \|\nabla u\|_{2}^2+\frac{t^{q(1+\delta_{q})}}{q} \|\nabla u\|_{q}^q-\frac{t^{ p\delta_{p}}}{p} {\|u\|}_p^p, \qquad\forall t\!>\!0.
$$
Differentiate $h(t)$ with respect to $t$, we obtain
$$h^{\prime}(t)=\frac{ t^{2} \|\nabla u\|_{2}^2+(1+\delta_{q})t^{q(1+\delta_{q})} \|\nabla u\|_{q}^q-\delta_{p}t^{ p\delta_{p}} {\|u\|}_p^p}{t}=\frac{P(u_t)}{t}.$$
Since $p\delta_{p}>2$, $p\delta_{p}>q(1+\delta_{q})$, $h^{\prime}(t)>0$ when $t>0$ is sufficiently small and $\lim_{t\to{ +\infty }}h^{\prime}(t)=-\infty$,  we derive that $h(t)$ has a unique maximum at some point $t_0 > 0$, see \cite{LLY}. Therefore, $h^{\prime}(t_0)=P(u_{t_0})/{t_0}=0$, which implies that $u_{t_0} \in \mathcal{P}_{c}$. In addition, we have  $I(u_{t_0} ) = \max_{\{t>0\}}I(u_t)$ and $P(u_{t_0})=0$. Next, we prove the last two statements. We first claim that $P(u)<0 \Rightarrow t_0< 1$. 
Just suppose that $t_0\geq1$, by using $h^{\prime}(t_0)=0$ and $P(u)<0$, we obtain the following contradiction
\begin{align*}
0&=t_{0}^{2-p\delta_{p}} \|\nabla u \|_{2}^2+(1+\delta_{q})t_{0}^{q(1+\delta_{q})-p\delta_{p}} \|\nabla u \|_{q}^q-\delta_{p} {\|u \|}_p^p  \\
&\leq \|\nabla u \|_{2}^2+(1+\delta_{q}) \|\nabla u \|_{q}^q-\delta_{p} {\|u \|}_p^p<0,
\end{align*}
so $t_0< 1$ is proved.
If $P(u)\!=\!0$, it is easy to show that neither $t_0>1$ nor $t_0<1$ could occur, thus $P(u)=0 \Rightarrow t_0= 1$. Next, we show that $t_0< 1  \Rightarrow P(u)<0$ and $t_0= 1  \Rightarrow P(u)=0$. If $t_0<1$, we see that
\begin{align*}
0&=t_{0}^{2-p\delta_{p}} \|\nabla u\|_{2}^2+(1+\delta_{q})t_{0}^{q(1+\delta_{q})-p\delta_{p}} \|\nabla u\|_{q}^q-\delta_{p} {\|u\|}_p^p  \\
&> \|\nabla u\|_{2}^2+(1+\delta_{q}) \|\nabla u\|_{q}^q-\delta_{p} {\|u\|}_p^p=P(u).
\end{align*}
Also, $t_0=1$ implies $P(u)=P(u_{t_0})=0$.
\end{proof}

\begin{lemma} \label{3.4pqlalplace} Let $c>0$ and and assume \eqref{h1.3a}, \eqref{h1.3b}. Then each minimizer of $I|_{\mathcal{P}_{c}}$
is a critical point of $I|_{S_{c}}$.
\end{lemma}

\begin{proof}
Suppose that $u$ is a minimizer of $I|_{\mathcal{P}_{c}}$, then $P(u)=\|\nabla u\|_{2}^2+(1+\delta_{q}) \|\nabla u\|_{q}^q-\delta_{p} {\|u\|}_p^p=0$. From \cite[Corollary 4.1.2]{Chan}, there exist two Lagrange multipliers $\lambda$ and $\mu$ such that
$$I'(u)-\lambda u-\mu P'(u)=0~~~~\mbox{in}~~~~X^{-1}.$$
That is to say, $u$ satisfies
$$
(2\mu-1)\Delta u+[q(1+\delta_{q})\mu-1]\Delta_q u+ (p\delta_{p}\mu-1){| u |^{p - 2}}u-\lambda u=0.$$
It is sufficient to prove that $\mu=0$.
In the same way as Lemma \ref{lem2.7}, we get
\begin{align} \label{pqeqA5.10}
(1-2\mu)\|\nabla u\|_{2}^2+[1-q(1+\delta_{q})\mu](1+\delta_{q}) \|\nabla u\|_{q}^q-(1-p\delta_{p}\mu)\delta_{p} {\|u\|}_p^p=0.
\end{align}
Recalling that $P(u)=0$, so \eqref{pqeqA5.10} can be reduced to
$$\mu\Big\{2\|\nabla u\|_{2}^2+q(1+\delta_{q})^2\|\nabla u\|_{q}^q-p\delta^2_{p} {\|u\|}_p^p\Big\}=0.$$
Using $P(u)=0$ again, we obtain
$$\mu\Big\{(p\delta_{p}-2)\|\nabla u\|_{2}^2+(1+\delta_{q})[p\delta_{p}-q(1+\delta_{q})]\|\nabla u\|_{q}^q\Big\}=0.$$
Since $p\delta_{p}>2$ and $p\delta_{p}>q(1+\delta_{q})$, we conclude that $\mu=0$.
\end{proof}

Lemma \ref{3.4pqlalplace} indicates that the restriction $P(u)=0$ in $\mathcal{P}_{c}$ is a natural constraint. To prove $\sigma(c)$ is attained, we show that $\sigma(c)$ is strictly decreasing with respect to $c$.

\begin{lemma} \label{3.5pqlalplace}
Let $c>0$ and assume \eqref{h1.3a}, \eqref{h1.3b}. If $c_2>c_1>0$, then $\sigma(c_2)<\sigma(c_1)$.
\end{lemma}

\begin{proof}
From Lemma \ref{pqlalpl3.1}, we deduce that $\sigma(c)>0$ for any $c>0$. By using Lemma \ref{3.2pqlalplace}, there exists a sequence $\{u_n\}\subset\mathcal{P}_{c_1}$ such that
$$\sigma(c_1) \leq I(u_{n}) = \max_{\{t>0\}}I(t^{\frac{N}{2}}u_n(tx))<\sigma(c_1)+\frac{1}{n}.$$


Furthermore, for each $u_n \in \mathcal{P}_{c_1}$, we have, like in \eqref{pq4.4},
\begin{align} \label{pq41.414}
\biggl(\frac{1}{2}-\frac{1}{p\delta_{p}}\biggr) \|\nabla u_{n}\|_{2}^2+\biggl(\frac{1}{q}-\frac{1+\delta_{q}}{p\delta_{p}}\biggr)\|\nabla u_{n}\|_{q}^q=I(u_{n})\leq\sigma(c_1)+1.
\end{align}
Consequently, $\{u_n\}$ is bounded in $X$ as $p\delta_{p}>2$ and $p\delta_{p}>q(1+\delta_{q})$. In a fashion similar to \eqref{pq4.3.1}, we also derive a positive lower bound for $\|\nabla u_n\|_{2}$, $\|\nabla u_n\|_{q}$ and $\|u_n\|_{p}$. Now, for $c_2>c_1>0$, we prove $\sigma(c_2)<\sigma(c_1)$ in two cases respectively.

\textbf{Case (i):} $N\!\geq\!2$, $\frac{2N(N+2)}{N^2+2N+4}\!<\!q\!<\!2$ and  $2(1\!+\!\frac{2}{N})\!<\!p\!<\!q^*$. Denote
$$\theta:=\biggl(\frac{c_2}{c_1}\biggr)^{\frac{2q}{q(N+2)-2N}}=\biggl(\frac{c_2}{c_1}\biggr)^{\frac{1}{1+\delta_q}}>1\quad \text{and}\quad v_n(x):=\theta^{\frac{q-N}{q}}u_n(\theta^{-1}x).$$
Direct computations give
\begin{align*}
&\|v_n\|_{2}^{2}=\theta^{N+2-\frac{2N}{q}}\|u_n\|_{2}^{2}=\theta^{2(1+\delta_q)}c_1^{2}=c_2^{2},\qquad\|\nabla v_n\|_{q}^{q}=\|\nabla u_n\|_{q}^{q},  \\
&\|\nabla v_n\|_{2}^{2}=\theta^{N-\frac{2N}{q}}\|\nabla u_n\|_{2}^{2}=\theta^{2\delta_q}\|\nabla u_n\|_{2}^{2},\qquad\|v_n\|_{p}^{p}=\theta^{N+p-\frac{Np}{q}}\|u_n\|_{p}^{p}.
\end{align*}
Using Lemma \ref{3.2pqlalplace} again, there exists $t_n > 0$ such that $t_n^{\frac{N}{2}}v_n(t_nx)\in\mathcal{P}_{c_2}$ and
$$I(t_n^{\frac{N}{2}}v_n(t_nx)) = \max_{t>0} I(t^{\frac{N}{2}}v_n(tx)).$$
Notice that $\{u_n\}$ is bounded, $P(t_n^{\frac{N}{2}}v_n(t_nx))=0$ and we could choose a positive constant $C>0$ such that $t_n>C>0$. It is easy to check that $2\delta_q<0$ and $N+p-\frac{Np}{q}>0$. Therefore, we have
\begin{align*}
\sigma(c_2)&\leq I(t_n^{\frac{N}{2}}v_n(t_nx))=\frac{t_n^{2}}{2}{\|\nabla v_n\|}_2^2+\frac{t_n^{q(1+\delta_q)}}{q}\|\nabla v_n\|_{q}^q-\frac{t_n^{p\delta_p}}{p}{\|v_n\|}_p^p\\
&=\frac{t_n^{2}}{2}\theta^{2\delta_q}{\|\nabla u_n\|}_2^2+\frac{t_n^{q(1+\delta_q)}}{q}\|\nabla u_n\|_{q}^q-\frac{t_n^{p\delta_p}}{p}\theta^{N+p-\frac{Np}{q}}{\|u_n\|}_p^p\\
&=I(t_n^{\frac{N}{2}}u_n(t_nx))+\frac{t_n^{2}}{2}\Big(\theta^{2\delta_q}-1\Big){\|\nabla u_n\|}_2^2+\frac{t_n^{p\delta_p}}{p}\Big(1-\theta^{N+p-\frac{Np}{q}}\Big) {\|u_n\|}_p^p\\
&\leq I(u_n)+\frac{t_n^{2}}{2}\Big(\theta^{2\delta_q}-1\Big){\|\nabla u_n\|}_2^2+\frac{t_n^{p\delta_p}}{p}\Big(1-\theta^{N+p-\frac{Np}{q}}\Big) {\|u_n\|}_p^p\\
&\leq\sigma(c_1)-C+\frac{1}{n}<\sigma(c_1)
\end{align*}
for $n$ sufficiently large.

\textbf{Case (ii):} $N\!\geq\!3$,  $2\!<\!q\!<\!\min\bigl\{ N,\frac{2N  ^2}{ {N}^{2}-4} \bigr\}$ and $q(1+\frac{2}{N})\!<\!p\!<\!2^*$. Let
$$w_n(x):=\biggl(\frac{c_2}{c_1}\biggr)^{1-\frac{N}{2}}u_n\biggl(\frac{c_1}{c_2}x\biggr),$$
then we deduce that
\begin{align*}
&\|w_n\|_{2}^{2}=\biggl(\frac{c_2}{c_1}\biggr)^{2}\|u_n\|_{2}^{2}
=c_2^2,\qquad\|\nabla w
_n\|_{q}^{q}=\biggl(\frac{c_2}{c_1}\biggr)^{N-\frac{Nq}{2}}\|\nabla u_n\|_{q}^{q}=\biggl(\frac{c_2}{c_1}\biggr)^{-q\delta_q}\|\nabla u_n\|_{q}^{q},  \\
&\|\nabla w_n\|_{2}^{2}=\|\nabla u_n\|_{2}^{2},\qquad\|w_n\|_{p}^{p}=\biggl(\frac{c_2}{c_1}\biggr)^{N+p-\frac{Np}{2}}\|u_n\|_{p}^{p}
=\biggl(\frac{c_2}{c_1}\biggr)^{p(1-\delta_p)}\|u_n\|_{p}^{p}.
\end{align*}
From Lemma \ref{3.2pqlalplace}, there exists $t_n > 0$ such that $t_n^{\frac{N}{2}}w_n(t_nx)\in\mathcal{P}_{c_2}$ and
$$I(t_n^{\frac{N}{2}}w_n(t_nx)) = \max_{t>0} I(t^{\frac{N}{2}}w_n(tx)).$$
As $\{u_n\}$ is bounded and $P(t_n^{\frac{N}{2}}w_n(t_nx))=0$, we could choose a positive constant $C>0$ such that $t_n>C>0$. Therefore, we have
\begin{align*}
\sigma(c_2)&\leq I(t_n^{\frac{N}{2}}w_n(t_nx))=\frac{t_n^{2}}{2}{\|\nabla w_n\|}_2^2+\frac{t_n^{q(1+\delta_q)}}{q}\|\nabla w_n\|_{q}^q-\frac{t_n^{p\delta_p}}{p}{\|w_n\|}_p^p\\
&=\frac{t_n^{2}}{2}{\|\nabla u_n\|}_2^2+\frac{t_n^{q(1+\delta_q)}}{q}\biggl(\frac{c_2}{c_1}\biggr)^{-q\delta_q}\|\nabla u_n\|_{q}^q-\frac{t_n^{p\delta_p}}{p}\biggl(\frac{c_2}{c_1}\biggr)^{p(1-\delta_p)}{\|u_n\|}_p^p\\
&=I(t_n^{\frac{N}{2}}u_n(t_nx))+\frac{t_n^{q(1+\delta_q)}}{q}
\biggl(\biggl(\frac{c_2}{c_1}\biggr)^{-q\delta_q}-1\biggr){\|\nabla u_n\|}_q^q+\frac{t_n^{p\delta_p}}{p}\biggl(1-\biggl(\frac{c_2}{c_1}\biggr)^{p(1-\delta_p)}\biggr) {\|u_n\|}_p^p\\
&\leq I(u_n)+\frac{t_n^{q(1+\delta_q)}}{q}
\biggl(\biggl(\frac{c_2}{c_1}\biggr)^{-q\delta_q}-1\biggr){\|\nabla u_n\|}_q^q+\frac{t_n^{p\delta_p}}{p}\biggl(1-\biggl(\frac{c_2}{c_1}\biggr)^{p(1-\delta_p)}\biggr) {\|u_n\|}_p^p\\
&\leq\sigma(c_1)-C+\frac{1}{n}<\sigma(c_1)
\end{align*}
for $n$ sufficiently large.

\end{proof}

Then, we prove that $\sigma(c)$ can be attained.
\begin{lemma} \label{5.5pqlalplace}
Let $c>0$ and assume \eqref{h1.3a}, \eqref{h1.3b}. Then $\sigma(c)\!=\!\inf _{u \in \mathcal{P}_{c}} I(u)$ is attained.
\end{lemma}

\begin{proof}
Let $\{u_n\}$ be a minimizing sequence for $\sigma(c)$. In a fashion similar to \eqref{pq41.414}, we see that $\{u_n\}$ is bounded in $X$. Then there exists $u\in X$ such that
\begin{align*}
u_{n}\rightharpoonup u~~~\mbox{in}~~~X,~~~~u_{n} \rightarrow u~~~\mbox{in}~~~L_{loc}^p({\R}^N),~~~~u_{n} \rightarrow u~~~\mbox{a.e. on}~~~{\R}^N
\end{align*}
up to a subsequence. We claim that $u\not\equiv0$. Otherwise, if $u\equiv0$, we learn from \cite[Lemma I.1]{LilP} that $u_{n}\rightarrow 0$ in $L^p({\R}^N)$, this along with $P(u_{n})\!=\! \|\nabla u_{n}\|_{2}^2+(1+\delta_{q})  \|\nabla u_{n}\|_{q}^q-\delta_{p}  {\|u_{n}\|}_p^p\!=\!0$ lead to $\|\nabla u_{n}\|_{2}\!\to\!0$ and $\|\nabla u_{n}\|_{q}\!\to\!0$. We obtain $\sigma(c)=0$, which contradicts Lemma \ref{pqlalpl3.1}, namely that $\sigma(c)>0$. Now let ${|u_{n}|}^*$ be the symmetric decreasing rearrangement of $|u_{n}|$ and denote $v_{n}={|u_{n}|}^*$. Similar to the proof of Theorem \ref{th1.1}, we have
$$\|\nabla v_{n}\|_2  \leq  \|\nabla{|u_{n}|} \|_2 \leq  \|\nabla{u_{n}} \|_2 ,~~~~~~\|\nabla v_{n} \|_q  \leq  \|\nabla{|u_{n}|} \|_q \leq  \|\nabla{u_{n}} \|_q ,~~~~~~\|v_{n}\|_p=\|u_{n}\|_p,$$
so we get $I\left(v_{n}\right)\!\leq\!I\left(|u_{n}|\right)\!\leq\!I\left(u_{n}\right)$ and $P\left(v_{n}\right)\!\leq\!P\left(|u_{n}|\right)\!\leq\!P\left(u_{n}\right)=0$. Utilizing Lemma \ref{3.2pqlalplace}, there exists $t_n \in (0, 1]$ such that
$$w_{n}(x):=t_n^{\frac{N}{2}}v_n(t_nx) \in \mathcal{P}_{c},~~~~P(w_{n})= 0.$$
Since $w_n \in \mathcal{P}_{c}$, $t_n \in (0, 1]$ and $u_n\in \mathcal{P}_{c}$, we have
$$\begin{aligned}
I(w_{n})&=\biggl(\frac{1}{2}-\frac{1}{p\delta_{p}}\biggr) \|\nabla w_{n}\|_{2}^2+\biggl(\frac{1}{q}-\frac{1+\delta_{q}}{p\delta_{p}}\biggr)\|\nabla w_{n}\|_{q}^q
\nonumber \\
&=\biggl(\frac{1}{2}-\frac{1}{p\delta_{p}}\biggr){t_n^2} \|\nabla v_{n}\|_{2}^2+\biggl(\frac{1}{q}-\frac{1+\delta_{q}}{p\delta_{p}}\biggr)t_n^{q(1+\delta_q)}\|\nabla v_{n}\|_{q}^q \nonumber \\
&\leq\biggl(\frac{1}{2}-\frac{1}{p\delta_{p}}\biggr) \|\nabla v_{n}\|_{2}^2+\biggl(\frac{1}{q}-\frac{1+\delta_{q}}{p\delta_{p}}\biggr)\|\nabla v_{n}\|_{q}^q \nonumber \\
&\leq\biggl(\frac{1}{2}-\frac{1}{p\delta_{p}}\biggr) \|\nabla u_{n}\|_{2}^2+\biggl(\frac{1}{q}-\frac{1+\delta_{q}}{p\delta_{p}}\biggr)\|\nabla u_{n}\|_{q}^q =I(u_{n}).
\end{aligned}$$
Thus, we obtain a new minimizing sequence $\{w_n\}$ for $\sigma(c)$. In addition, we can prove that, $\{w_n\}$ is bounded in $X$ and there exists $w\in X_{r}\!:=\!\{u(x) \!\in\! X: u(x)\!=\!u(|x|)\}$ such that
\begin{align}  \label{pq41.417}
w_{n}\rightharpoonup w\not\equiv0~~~\mbox{in}~~~X_r,~~~~w_{n} \rightarrow w~~~\mbox{in}~~~L^p_{loc}({\R}^N),~~~~w_{n} \rightarrow w~~~\mbox{a.e. on}~~~{\R}^N
\end{align}
up to a subsequence. We claim that $\|w\|_{2}=c$. Otherwise, if $\|w\|_{2}=c_1\in(0,c)$, then Lemma \ref{3.5pqlalplace} indicates $\sigma(c)<\sigma(c_1)$. By using \eqref{pq41.417}, we have $P\left(w\right)\!\leq\!\lim_{n\to\infty}P\left(w_{n}\right)=0$. From Lemma \ref{3.2pqlalplace}, there exists $\tau_0 \in (0, 1]$ such that
$$\tau_0^{\frac{N}{2}}w(\tau_0x) \in \mathcal{P}_{c_1},~~~~P(\tau_0^{\frac{N}{2}}w(\tau_0x))= 0.$$
However, we have a contradiction from
\begin{align} \label{pq41.418}
\sigma(c)<\sigma(c_1)&\leq I(\tau_0^{\frac{N}{2}}w(\tau_0x))=\biggl(\frac{1}{2}-\frac{1}{p\delta_{p}}\biggr){\tau_0^2} \|\nabla w\|_{2}^2+\biggl(\frac{1}{q}-\frac{1+\delta_{q}}{p\delta_{p}}\biggr)\tau_0^{q(1+\delta_q)}\|\nabla w\|_{q}^q \nonumber \\
&\leq\biggl(\frac{1}{2}-\frac{1}{p\delta_{p}}\biggr) \|\nabla w\|_{2}^2+\biggl(\frac{1}{q}-\frac{1+\delta_{q}}{p\delta_{p}}\biggr)\|\nabla w\|_{q}^q \nonumber \\
&\leq\biggl(\frac{1}{2}-\frac{1}{p\delta_{p}}\biggr) \lim_{n\to\infty} \|\nabla w_n\|_{2}^2+\biggl(\frac{1}{q}-\frac{1+\delta_{q}}{p\delta_{p}}\biggr)\lim_{n\to\infty}\|\nabla w_n\|_{q}^q \nonumber \\
&= \lim_{n\to\infty}I\left(w_{n}\right)=\sigma(c).
\end{align}
So, it must be that $c_1=c$ and $\tau_0=1$. That is to say, $\|w\|_{2}=c$ and $\sigma(c)=I\left(w\right)$. From \eqref{pq41.418} and $\|w\|_{2}=c$, we also obtain $w_{n}\rightarrow w$ in $X_r$ and $w$ is a minimizer for $\sigma(c)$.

\end{proof}

\begin{proof}[\textbf {Proof of Theorem \ref{th1.3}.} ] From Lemma \ref{5.5pqlalplace}, we see that $\sigma(c)$ is attained by some $u_c \!\in\! S_c $. If $v \in S_c$ and $(I|_{S_c})'(v)=0$, we deduce from Lemma \ref{lem2.7} that $v \in \mathcal{P}_{c}$. Therefore, we have $I(v)\geq\sigma(c)=I(u_c)$. This implies that $u_c$ is a ground state of \eqref{eq1.1}. Next, the Lagrange multipliers rule indicates the existence of some $\lambda_c\in \mathbb{R}$ such that
$$\int_{{\R^N}}[\nabla u_c\nabla\varphi  \!+\! {| \nabla u_c |^{q - 2}}\nabla u_c \nabla{\varphi}  \!-\! \left|u_c\right|^{p\!-\!2} u_c {\varphi}  \!-\!\lambda_c  u_c{\varphi} ]dx=0,~~~~\forall \varphi \in X.$$
That is, $(\lambda_c,u_c)$ satisfies \eqref{eq1.1}. Under the assumptions \eqref{h1.3a}, \eqref{h1.3b}, we can check that
$$
1-\frac{1}{\delta_{p}}<0,\qquad 1-\frac{ 1+\delta_{q} }{\delta_{p}}<0
$$
and
$$
\frac{p(1-\delta_{p})}{p\delta_{p}-2}>0,\qquad\frac{p(1-\nu_{p,q})}{p\nu_{p,q}-q}>0,
\qquad\frac{1}{2}-\frac{1}{p\delta_{p}}>0,\qquad 1-\frac{q(1+\delta_{q})}{p\delta_{p}}>0.
$$
So we get
\begin{equation}\label{5.19bi}
\lambda_cc^2=\|\nabla u_c\|_{2}^2+\|\nabla u_c\|_{q}^q-{\|u_c\|}_p^p=\biggl(1-\frac{1}{\delta_{p}}\biggr)\|\nabla u_c\|_{2}^2+\biggl(1-\frac{ 1+\delta_{q} }{\delta_{p}}\biggr)\|\nabla u_c\|_{q}^q<0,
\end{equation}
which gives $\lambda_c<0$.
In the same way as in \eqref{pq4.3.1}, we have
\begin{align*}
\left\|\nabla u_c \right\|_{2}\geq\biggl[\frac{1}{\delta_{p} \mathcal{C}^p_{N,p}}\biggr]^{ \frac{1}{p\delta_{p}-2} }c^{-\frac{p(1-\delta_{p})}{p\delta_{p}-2}}~~\mbox{and}~~
\left\|\nabla u_c \right\|_{q}\geq\biggl[\frac{1+\delta_{q}}{\delta_{p} \mathcal{K}^p_{N,p}}\biggr]^{ \frac{1}{\nu_{p,q}-q} }c^{-\frac{p(1-\nu_{p,q})}{p\nu_{p,q}-q}}.
\end{align*}
By using $u_c \in \mathcal{P}_{c}$, we can rewrite $I(u_c)$ in the form
\begin{align*} 
\sigma(c)&=I(u_c)=\frac{1}{2} \|\nabla u_c\|_{2}^2+\frac{1}{q} \|\nabla u_c\|_{q}^q-\frac{1}{p} {\|u_c\|}_p^p  \nonumber \\
&=\biggl(\frac{1}{2}-\frac{1}{p\delta_{p}}\biggr) \|\nabla u_c\|_{2}^2+\biggl(\frac{1}{q}-\frac{1+\delta_{q}}{p\delta_{p}}\biggr)\|\nabla u_c\|_{q}^q \nonumber \\
&\geq \biggl(\frac{1}{2}-\frac{1}{p\delta_{p}}\biggr) \biggl[\frac{1}{\delta_{p} \mathcal{C}^p_{N,p}}\biggr]^{ \frac{2}{p\delta_{p}-2} }c^{-\frac{2p(1-\delta_{p})}{p\delta_{p}-2}} +\frac{1}{q}
\biggl(1-\frac{q(1+\delta_{q})}{p\delta_{p}}\biggr)  \Big[\frac{1+\delta_{q}}{\delta_{p} \mathcal{K}^p_{N,p}}\Big]^{ \frac{q}{p\nu_{p,q}-q} }c^{-\frac{qp(1-\nu_{p,q})}{p\nu_{p,q}-q}} .
\end{align*}
Therefore, we have $\sigma(c) \gtrsim c^{-\frac{2p(1-\delta_{p})}{p\delta_{p}-2}}+c^{-\frac{qp(1-\nu_{p,q})}{p\nu_{p,q}-q}}$. This implies that $\sigma(c) \to +\infty$ as $c \to 0^{+}$. Similarly, from \eqref{5.19bi}, we have $\lambda_c \lesssim -\bigl(c^{-\frac{2(p-2)}{p\delta_p-2}}+c^{ -\frac{2q(p-2)}{Np-Nq-2q}}\bigr)$. Thus, we see that $\lambda_c \to -\infty$ as $c \to 0^{+}$.
\end{proof}

\subsection{Existence of infinitely many radial solutions} \label{subsec5.2}

In this subsection, we study the existence of infinitely many radial solutions.

Recall that $X_r\!=\!\{u(x) \!\in\! H^{1}({\mathbb{R}^N}) \!\cap\! D^{1,q}({\mathbb{R}^N}): u(x)\!=\!u(|x|)\}$. Let $\{ {V_n}\}  \subset X_r$ be a strictly increasing sequence of finite-dimensional linear subspace in $X_r$ such that ${\bigcup\limits_n {{V_n}} }$ is dense in $X_r$. Denote $V_n^ \bot$ the orthogonal space of $V_n$ in $X_r$, then we have
\begin{lemma} \label{LemMa5.1}
Assume \eqref{h1.3a}, \eqref{h1.3b}. Then there holds
\[{\mu _n}: = \mathop {\inf }\limits_{u \in V_{n - 1}^ \bot } \frac{ (\|\nabla u\|_{2}+\|\nabla u\|_{q}+{\|u\|}_2)^2 }{ {\|u\|}_p^2  } = \mathop {\inf }\limits_{u \in V_{n - 1}^ \bot } \frac{{{{\left\| u \right\|}_X^2}}}{{\left\| u \right\|_p^2}} \to +\infty \text{ as } n \to +\infty .\]
\end{lemma}

\begin{proof}
It is similar to that of \cite[Lemma 2.1]{TSdV}.
\end{proof}

From now on, let $c>0$ be fixed and
$$S_{r,c}:=\{u\in S_{c},~~u(x)=u(|x|)\}.$$
For any $n \in {\mathbb{N}^+ }$ and $n\geq 2$, we define
\[
{\rho _n}: = \min\biggl\{ \biggl[\frac{ p{\mu}_n^{\frac{p}{2}} }{ 2^{p+1}L }\biggr]^{\frac{1}{p-2}}, \biggl[\frac{ p{\mu}_n^{\frac{p}{2}} }{ 2^{p} qL }\biggr]^{\frac{1}{p-q}} \biggr\} \text{ with }L = \mathop {\max }\limits_{x > 0} \frac{(x+c)^p}{x^p+c^p},
\]

\begin{equation}\label {bu6.1}
\begin{array}{rl}
\displaystyle
{B_n}: = \{ u \in V_{n - 1}^ \bot  \cap { S_{r,c} }:\|\nabla u\|_{2}+\|\nabla u\|_{q}= {\rho _n}\},
\end{array}
\end{equation}
and

\begin{equation}\label {beta6.2}
\begin{array}{rl}
\displaystyle
{\beta_n}: = \mathop {\inf }\limits_{u \in {B_n}} I(u).
\end{array}
\end{equation}
In the following, we use a subtle analysis to prove that ${\beta_n} \to +\infty \text{ as }n \to +\infty$.
\begin{lemma} \label{6.2.5.2}
Assume \eqref{h1.3a}, \eqref{h1.3b}, then ${\beta_n} \to +\infty \text{ as }n \to +\infty$.
\end{lemma}
\begin{proof}
For any fixed $u \in {B_n}$, we have $\|\nabla u\|_{2}+\|\nabla u\|_{q}= {\rho _n}$ and
$$\|\nabla u\|_{2}\leq{\rho _n}\leq\biggl[\frac{ p{\mu}_n^{\frac{p}{2}} }{ 2^{p+1} L }\biggr]^{\frac{1}{p-2}},\qquad \|\nabla u\|_{q}\leq{\rho _n}\leq\biggl[\frac{ p{\mu}_n^{\frac{p}{2}} }{ 2^{p} qL }\biggr]^{\frac{1}{p-q}}.$$
Consequently, $(v+w)^{\gamma}\le 2^{\gamma-1}(v^\gamma+w^\gamma)$
for  $\gamma=p\ge1$, we obtain
$$\begin{aligned}
I(u) =&\frac{1}{2} \|\nabla u\|_{2}^2+\frac{1}{q} \|\nabla u\|_{q}^q-\frac{1}{p} {\|u\|}_p^p    \nonumber \\
\ge& \frac{1}{2} \|\nabla u\|_{2}^2+\frac{1}{q} \|\nabla u\|_{q}^q- \frac{1}{{ p{\mu} _n^{\frac{p}{2}} }}{(\|\nabla u\|_{2}+\|\nabla u\|_{q}+c)^{p} } \nonumber \\
\ge& \frac{1}{2} \|\nabla u\|_{2}^2+\frac{1}{q} \|\nabla u\|_{q}^q- \frac{L}{{ p{\mu} _n^{\frac{p}{2}} }}{  \Big[(\|\nabla u\|_{2}+\|\nabla u\|_{q})^{p}+c^{p}\Big] } \nonumber \\
\ge& \frac{1}{2} \|\nabla u\|_{2}^2+\frac{1}{q} \|\nabla u\|_{q}^q- \frac{ 2^{p-1} L }{{ p{\mu}_n^{\frac{p}{2}} }}{  (\|\nabla u\|_{2}^{p}+\|\nabla u\|_{q}^{p}) }- \frac{Lc^{p}}{{ p{\mu}_n^{\frac{p}{2}} }}    \nonumber  \\
=& \frac{1}{2} \|\nabla u\|_{2}^2- \frac{ 2^{p-1} L }{{ p{\mu}_n^{\frac{p}{2}} }}\|\nabla u\|_{2}^{2} \|\nabla u\|_{2}^{p-2}+\frac{1}{q} \|\nabla u\|_{q}^q- \frac{ 2^{p-1} L }{{ p{\mu}_n^{\frac{p}{2}} }}\|\nabla u\|_{q}^{q}\|\nabla u\|_{q}^{p-q}- \frac{Lc^{p}}{{ p{\mu}_n^{\frac{p}{2}} }} \nonumber  \\
\geq &\frac{1}{2} \|\nabla u\|_{2}^2- \frac{ 2^{p-1} L }{{ p{\mu}_n^{\frac{p}{2}} }}\|\nabla u\|_{2}^{2} \biggl[\frac{ p{\mu}_n^{\frac{p}{2}} }{ 2^{p+1} L }\biggr]   +\frac{1}{q} \|\nabla u\|_{q}^q- \frac{ 2^{p-1} L }{{ p{\mu}_n^{\frac{p}{2}} }}\|\nabla u\|_{q}^{q}\biggl[\frac{ p{\mu}_n^{\frac{p}{2}} }{ 2^{p} qL }\biggr]- \frac{Lc^{p}}{{ p{\mu} _n^{\frac{p}{2}} }} \nonumber   \\
\geq & \frac{1}{4} \|\nabla u\|_{2}^2+\frac{1}{2q} \|\nabla u\|_{q}^q-\frac{Lc^{p}}{{ p{\mu}_n^{\frac{p}{2}} }}.
\end{aligned}$$
By using Lemma \ref{LemMa5.1}, we have ${\mu_n} \to +\infty $ as $n\rightarrow+\infty$. Since $p>2$ and $p>q$, we derive from the definition of ${\rho_n}$ that ${\rho_n}\to +\infty $ as $n\rightarrow+\infty$. Recall that, we have $\|\nabla u\|_{2}+\|\nabla u\|_{q}= {\rho _n}$ for any $u \in {B_n}$, so it must be that at least one of $\|\nabla u\|_{2}$ and $\|\nabla u\|_{q}$ goes to $+\infty $ as $n\rightarrow+\infty$. All the above facts indicate that, for any fixed $u \in {B_n}$, it holds that
\begin{align*}
I(u) \geq \frac{1}{4} \|\nabla u\|_{2}^2+\frac{1}{2q} \|\nabla u\|_{q}^q-\frac{Lc^{p}}{{ p{\mu}_n^{\frac{p}{2}} }} \to+\infty
\end{align*}
as $n\rightarrow+\infty$. At last, take infimum with respect to $u \in {B_n}$, we have ${\beta_n}: = \mathop {\inf }\limits_{u \in {B_n}} I(u)\to+\infty$.
\end{proof}

Next, we define a map
\begin{equation}\label{5.4-uthte}
\kappa: X_r \times \R \to X_r,~~~~(u,\theta ) \mapsto  \kappa(u,\theta )=(u\star\theta): = {e^{{\textstyle{\frac{N}{2}}}\theta }}u({e^\theta }x).
\end{equation}
For any fixed $u \in S_{r,c}$, we have $u\star\theta\in S_{r,c}$ for all $\theta \in \R$ and that
\begin{equation*}
I(u\star\theta)=I({e^{{\textstyle{\frac{N}{2}}}\theta }}u({e^\theta }x))=\frac{e^{2\theta}}{2} \|\nabla u\|_{2}^2+\frac{e^{q(1+\delta_{q})\theta}}{q} \|\nabla u\|_{q}^q-\frac{e^{ p\delta_{p}\theta}}{p} {\|u\|}_p^p.
\end{equation*}
From \eqref{cond3}, we have $p\delta_{p}>2$ and $p\delta_{p}>q(1+\delta_{q})$, so that
$$
\begin{array}{rl}
\displaystyle
\left\{ {\begin{array}{*{20}{c}}
   {\|\nabla (u\star\theta) \|_{2}+\|\nabla (u\star\theta) \|_{q}\to 0^+,\text{ }I(u\star\theta) \to 0^+,\text{ as }\theta  \to  - \infty }  \hfill \\
   {\|\nabla (u\star\theta) \|_{2}+\|\nabla (u\star\theta) \|_{q} \to  + \infty,\text{ }I(u\star\theta) \to  - \infty ,\text{ as }\theta  \to  + \infty }.  \hfill \\
\end{array}} \right.
\end{array}$$
Hence, there exists $\theta_n>0$ sufficiently large such that
\begin{equation}\label {6.1}
\begin{array}{rl}
\displaystyle
\left\{ {\begin{array}{*{20}{c}}
   {  \|\nabla u\star(-{\theta _n}) \|_{2}+\|\nabla  u\star(-{\theta _n})  \|_{q} < {\rho _n}<\|\nabla (u\star {\theta_n}) \|_{2}+\|\nabla (u\star {\theta_n}) \|_{q},  }  \\
   {   \max\{I( u\star(-{\theta _n}) ),~~I( u\star {\theta_n} )\} < {\beta_n}}.  \\
\end{array}} \right.
\end{array}
\end{equation}
Let
\begin{align*}
{\Gamma _n}\!:=\! \Big\{ \gamma\!\in\! C\big([0,1] \!\times\! ( S_{r,c} \!\cap\! {V_n}),S_{r,c}\big)~~|~~\gamma\text{ is odd in }u,
\gamma(0,u)\!=\!u\star(-{\theta _n}), \gamma(1,u)\!=\!u\star {\theta_n} \Big\},
\end{align*}
then ${\Gamma _n}\not=\emptyset$ as $u\star(2t-1){\theta _n}\in{\Gamma _n}$.
So we can introduce a sequence of min-max values
\begin{align}\label {sigma6.1}
{\sigma_n}(c): = \mathop {\inf }\limits_{\gamma \in {\Gamma _n}} \mathop {\max }\limits_{0 \le t \le 1,u \in {S_{r,c}} \cap {V_n}} I(\gamma(t,u)).
\end{align}

The following lemma has been proved in \cite{TlYl} via the Borsuk-Ulam theorem.
\begin{lemma} \label{6.1..1}(\cite[Lemma 3.2]{TlYl})
Let ${L_1}$ and $L$ be two finite dimensional normed vector spaces such that ${L_1}\subseteq L$ and ${L_1}\not=L$. For $\alpha\in \R$ and $u\!\in\! S\!=\!\{u\!\in\! L:~~\|u\|_2\!=\!c\!>\!0\}$, if $\eta=(\eta_1,\eta_2)\in C \big([0,1]\times S, \R \times L_1\big)$ satisfies   $$\eta_1(t,u)=\eta_1(t,-u),~~\eta_2(t,u)=-\eta_2(t,-u),~~\eta_1(0,u)<\alpha<\eta_1(1,u), $$
then there exists $(t,u)\!\in\![0,1]\!\times\! S$ such that $\eta(t,u)\!=\!(\alpha,0)$.
\end{lemma}

We will use Lemma \ref{6.1..1} to prove the key intersection Lemma.
\begin{lemma} \label{6.1.inter}
For each $\gamma\in\Gamma_n $, there exists $(t,u)\in [0,1] \times ({S_{r,c}} \cap {V_n})$ such that $\gamma(t,u) \in {B_n}$ with ${B_n}$ defined in \eqref{bu6.1}. That is, for each $\gamma\in\Gamma_n $, we have
$$\gamma\big(  [0,1] \times ({S_{r,c}} \cap {V_n}) \big)\cap{B_n}\not=\emptyset.$$
\end{lemma}
\begin{proof}
Take $L=V_n$ and $L_1=V_{n-1}$ endowed with the $L^2(\R^N)$ norm, $S={S_{r,c}} \cap {V_n}$ and $\alpha=\rho_n$ in Lemma \ref{6.1..1}. Denote $P_{n-1}$ the orthogonal projection from $X_r$ to $V_{n-1}$ and let
\begin{align*}
 h_n : S_{r,c} &\to \R \times V_{n-1},\\
       u &\mapsto (\|\nabla u\|_{2}+\|\nabla u\|_{q},P_{n-1}u).
\end{align*}
For $\gamma\in\Gamma_n $, define
\begin{align*}
\eta=(\eta_1,\eta_2)=h_n\circ \gamma:~~~~[0,1]\times ({S_{r,c}} \cap {V_n}) &\to \R \times V_{n-1},\\
(t,u) &\mapsto  h_n\big(\gamma(t,u)\big).
\end{align*}
It is obvious that
$$
h_n(u)=(\rho_n,0) \Leftrightarrow u\in{B_n}\!:=\!\{ v \in V_{n - 1}^ \bot  \!\cap\! { S_{r,c} }:\|\nabla v\|_{2}\!+\!\|\nabla v\|_{q}\!=\!{\rho _n}\}.
$$
We can verify, by \eqref{6.1}, that all the conditions of Lemma \ref{6.1..1} are satisfied. Thus, there exists $(t,u)\in [0,1] \times ({S_{r,c}} \cap {V_n})$ such that $\eta(t,u)=(\rho_n,0)$. By the fact that $\eta(t,u)=h_n\big(\gamma(t,u)\big)$, we obtain $\gamma(t,u)\in B_n$.
\end{proof}

\begin{lemma} \label{6.1.sig-beta}
Let ${\beta_n}$ and ${\sigma_n}(c)$ be defined in \eqref{beta6.2} and \eqref{sigma6.1}, then ${\sigma_n}(c)\ge {\beta_n}$.
\end{lemma}
\begin{proof}
From Lemma \ref{6.1.inter},  we see that $\gamma\big(  [0,1] \times ({S_{r,c}} \cap {V_n}) \big)\cap{B_n}\not=\emptyset$ for each $\gamma\in\Gamma_n $. Therefore, we get
$$
{\sigma_n}(c): = \mathop {\inf }\limits_{\gamma \in {\Gamma _n}} \mathop {\max }\limits_{0 \le t \le 1,u \in {S_{r,c}} \cap {V_n}} I(\gamma(t,u))\geq {\beta_n}: = \mathop {\inf }\limits_{u \in {B_n}} I(u).
$$
\end{proof}

According to Lemma \ref{6.2.5.2} and Lemma \ref{6.1.sig-beta}, there exists $n_0 \in \mathbb{N}$ such that
\begin{equation}\label{snbn}
{\sigma_n}(c)\ge {\beta_n}>1~~~~\text{for}~~~~n\geq n_0.
\end{equation}
For any fixed $n\geq n_0$, we shall prove that the sequence $\{ {\sigma_n}(c) \} (n\geq n_0)$  is indeed a sequence of critical values for $I|_{S_{r,c}}$. To this end, we will construct a bounded Palais-Smale sequence at level ${\sigma_n}(c)$. First, we introduce the stretched functional
\[
\widetilde I:S_{r,c} \times \R \to \R,\text{ } (u,\theta ) \mapsto I( u\star\theta )
\]
and a new min-max value
\begin{equation}\label {sigma6.1.3}
{\widetilde\sigma _n}(c): = \mathop {\inf }\limits_{\widetilde \gamma \in {\widetilde \Gamma _n}} \mathop {\max }\limits_{0 \le t \le 1,u \in {S_{r,c}  \cap {V_n}}} \widetilde I(\widetilde \gamma(t,u)),
\end{equation}
where $u\star\theta$ is given in \eqref{5.4-uthte} and the min-max class is given by
$$
{\widetilde \Gamma _n}: = \Bigl\{ \widetilde \gamma\in C\big([0,1] \times (S_{r,c}\cap {V_n}), S_{r,c} \times \R\big)\text{ }|\text{ }\widetilde \gamma \text{ is odd in }u\text{ and }\kappa  \circ \widetilde \gamma \in {\Gamma _n}\Bigr\}.
$$

Clearly, for any $\gamma \in {\Gamma _n}$, $\widetilde \gamma:=(\gamma,0) \in {\widetilde \Gamma _n}$.

\begin{lemma} \label{6.1.2}
Let ${\sigma_n}(c)$ and ${\widetilde\sigma _n}(c)$ be defined in \eqref{sigma6.1} and \eqref{sigma6.1.3} respectively, then ${\widetilde\sigma _n}(c)={\sigma_n}(c)$.
\end{lemma}

\begin{proof} 
Since the maps
$$
\varphi :{\Gamma _n} \to {\widetilde\Gamma _n},\text{ } \gamma \mapsto \varphi (\gamma): = (\gamma,0)~~~~\text{and}~~~~\psi :{\widetilde\Gamma _n} \to {\Gamma _n},\text{ } \widetilde \gamma \mapsto \psi (\widetilde \gamma): = \kappa  \circ \widetilde \gamma
$$
satisfy $\widetilde I(\varphi (\gamma)) = I(\gamma)$ and $I(\psi (\widetilde \gamma)) = \widetilde I(\widetilde \gamma)$, we immediately obtain ${\widetilde\sigma _n}(c)={\sigma_n}(c)$.
\end{proof}

Denote $E:=X_r\times \R$ endowed with the norm $\left\|  \cdot  \right\|_E= {\left\|  \nabla \cdot  \right\|_2} +{\left\|  \nabla \cdot  \right\|_q}+ {\left\| \cdot  \right\|_2}+\left|  \cdot  \right|_\R$ and ${E^{-1}}$ its dual space. From Ekeland's variational principle, we get the following lemma.

\begin{lemma} \label{5.12}
Let $\varepsilon>0$. If ${\widetilde \gamma} \in {\widetilde\Gamma _n}$ satisfies \[\mathop {\max }\limits_{0 \le t \le 1,u \in { S_{r,c} } \cap {V_n}} \widetilde I({\widetilde \gamma}(t,u)) \le {\widetilde\sigma_n}(c) + \varepsilon, \]
then there exists a pair of $({u},{\theta}) \in S_{r,c} \times \R$ such that:

(i) $\widetilde I({u},{\theta}) \in [{\widetilde \sigma_n}(c) - \varepsilon ,{\widetilde  \sigma_n}(c) + \varepsilon]$;

(ii) $\mathop {\min }\limits_{0 \le t \le 1,u \in S_{r,c} \cap {V_n}} { \| {({u},{\theta}) - {{\widetilde \gamma}}(t,u)}  \|_E} \le \sqrt \varepsilon$;

(iii) For all $(\phi,s) \in {\widetilde T_{({u},{\theta})}}: = \{ (\phi,s) \in E,{ \langle {{u},{\phi}}  \rangle _{{L^2}}} = 0\}$, it holds that
$${ \| { ( {{{\widetilde I}}\left| {_{ S_{r,c} \times \R}} \right.}  )^{'}}({u},{\theta})  \|_{{E^{-1}  }}} \le 2\sqrt \varepsilon
,~~~~i.e.,~~~~  | {{{ \langle {{{\widetilde I}^{'}}({u},{\theta}),(\phi,s)}  \rangle }_{{E^{-1} } \times E}}} | \le 2\sqrt \varepsilon  { \| (\phi,s)  \|_E}.$$
\end{lemma}
\begin{proof}
The proof is similar to that of \cite[Lemma 2.3]{lJaE}.
\end{proof}

\begin{proposition} \label{5.7-ps}
Assume \eqref{h1.3a}, \eqref{h1.3b}. Then, for any fixed $c>0$ and $n\geq n_0$, with $n_0$ defined above in order to have \eqref{snbn}, there exists a sequence $\{ {v_k^{(n)}}\}\subset S_{r,c}$ such that
$$I( {v_k^{(n)} }) \to { \sigma_n}(c),~~~~\big( {{\left. {{I}} \right|}_{ S_{r,c} }} \big)^{'}({v_k^{(n)}}) \to 0,~~~~{P({v_k^{(n)}}) \to 0}~~~~\mbox{as}~~~~k \to +\infty.$$
\end{proposition}

\begin{proof}
For each $k \in {\mathbb{N}^ + }$, there exists a ${\gamma_k} \in {\Gamma _n}$ such that
 \[
 \mathop {\max }\limits_{0 \le t \le 1,u \in S_{r,c} \cap {V_n}} I({\gamma_k}(t,u)) \le \sigma_n(c) + \frac{1}{k}.
 \]
By using $\widetilde\sigma_n(c) = \sigma_n(c)$ and ${\widetilde \gamma_k} = ({\gamma_k},0) \in {\widetilde\Gamma _n}$, we have
 \[
 \mathop {\max }\limits_{0 \le t \le 1,u \in  S_{r,c} \cap {V_n}} \widetilde I({\widetilde \gamma_k}(t,u)) \le \widetilde\sigma_n(c) + \frac{1}{k}
 \]
Taking $\varepsilon=\frac{1}{k}$ and ${\widetilde \gamma} = {\widetilde \gamma_k} \in {\widetilde\Gamma _n}$ in Lemma \ref{5.12}, we obtain a sequence $\{ ({u_k^{(n)}},{\theta_k^{(n)}})\}  \subset S_{r,c} \times \R$ such that

$(a)~~\widetilde I({u_k^{(n)}},{\theta_k^{(n)}}) \in [ \sigma_n(c)-\frac{1}{k}, \sigma_n(c)+ \frac{1}{k}]$;

$(b)~~\mathop {\min }\limits_{0 \le t \le 1,u \in S_{r,c} \cap {V_n}} { \| {({u_k^{(n)}},{\theta _k^{(n)}}) - ({{\gamma}_k}(t,u)},0) \|_E} \le  \frac{1}{\sqrt k}$;

$(c)$~~For all $(\phi,s) \in {\widetilde T_{({u_k^{(n)}},{\theta _k^{(n)}})}}: = \{ (\phi,s) \in E,\text{ }{ \langle {{u_k^{(n)}},{\phi}} \rangle _{{L^2}}} = 0\}$, it holds that
$$
{ \| { ( {{{\widetilde I}}\left| {_{  S_{r,c} \times \R}} \right.} )^{'}}({u_k^{(n)}},{\theta _k^{(n)}})  \|_{{E^{-1}  }}} \le \frac{2}{\sqrt k},~~~~i.e.,~~~~  | {{{ \langle {{{\widetilde I}^{'}}({u_k^{(n)}},{\theta _k^{(n)}}),(\phi,s)}  \rangle }_{{E^{-1}  } \times E}}}  | \le \frac{2}{\sqrt  k}  { \| (\phi,s) \|_E}.
$$
Denote
$$
{v_k^{(n)}} = \kappa ({u_k^{(n)}},{\theta _k^{(n)}})={u_k^{(n)}}\star{\theta _k^{(n)}}.
$$
Utilizing $(a)$ and ${I({v_k^{(n)}}) \!=\! I(\kappa ({u_k^{(n)}},{\theta _k^{(n)}})) \!=\! \widetilde I({u_k^{(n)}},{\theta _k^{(n)}})}$, we have $I({v_k^{(n)}}) \!\to\! \sigma_n(c)$ as $k \!\to\! \infty $. Next, we use $(c)$ to prove $P({v_k^{(n)}}) \to 0$ as $k \to \infty$. In fact, for all $(\phi,s) \in {\widetilde T_{({u_k^{(n)}},{\theta _k^{(n)}})}} = \{ (\phi,s) \in E,\text{ }{ \langle {{u_k^{(n)}},{\phi}}  \rangle _{{L^2}}} = 0\}$, we deduce from $(c)$ that
\begin{align*}
 & \langle {{{\widetilde I}^{'}}({u_k^{(n)}},{\theta _k^{(n)}}), (\phi,s)}  \rangle=\frac{d}{dt} I\Bigl((u_k^{(n)}+t\phi)\star(\theta_k^{(n)}+ts)\Bigr)\biggl|_{t=0} \\
=&s\big[ {e^{ 2{\theta _k^{(n)}} }}  {{\| {\nabla {u_k^{(n)}} } \|}_2^2}+  (1+\delta_q){e^{ q(1+\delta_q){\theta _k^{(n)}} }}{{\| {\nabla {u_k^{(n)}} } \|}_q^q}  - \delta_p{e^{ p \delta_p {\theta _k^{(n)}} }}{\| { {u_k^{(n)}} } \|}_p^p \big]+{e^{ 2{\theta _k^{(n)}} }}\int_{{\R^N}} {\nabla {u_k^{(n)}} \nabla \phi } dx  \\
&+ {e^{ q(1+\delta_q){\theta _k^{(n)}} }}\int_{{\R^N}} { | {\nabla {u_k^{(n)}} } |}^{q-2} {\nabla {u_k^{(n)}} \nabla \phi }dx-{e^{ p \delta_p {\theta _k^{(n)}} }}\int_{{\R^N}} {{{ | {u_k^{(n)}}  |}^{p-2}} {u_k^{(n)}} \phi }dx=o_k(1) { \| (\phi,s)  \|_E}.
\end{align*}
In particular, take $(\phi,s)=(0,1)$ in the above equality, we have
\begin{align*}
o_k(1)&= \langle {{{\widetilde I}^{'}}({u_k^{(n)}},{\theta _k^{(n)}}),(0,1)} \rangle    \\
&={e^{ 2{\theta _k^{(n)}} }}  {{\| {\nabla {u_k^{(n)}} } \|}_2^2}+  (1+\delta_q){e^{ q(1+\delta_q){\theta _k^{(n)}} }}{{\| {\nabla {u_k^{(n)}} } \|}_q^q}  - \delta_p{e^{ p \delta_p {\theta _k^{(n)} } }}{\| { {u_k^{(n)} } } \|}_p^p   \\
&={{\| {\nabla {v_k^{(n)}} } \|}_2^2}+  (1+\delta_q){{\| {\nabla {v_k^{(n)}} } \|}_q^q}  - \delta_p  {\| { {v_k^{(n)}} } \|}_p^p=P({v_k^{(n)}}).
\end{align*}
Thus, we have $P({v_k^{(n)}}) \to 0$ as $k \to \infty$ that is the claim. At last, we claim that
$$
{({\left. {{I}} \right|}_{{S_{r,c}}})^{'}({v_k^{(n)}}) \to 0 \text{ as }k \to \infty }.
$$
It is sufficient to prove that
$$
 | { \langle {{I^{'}}({v_k^{(n)}}),\phi } \rangle } | \le \frac{4}{{\sqrt k }}{ \| \phi  \|_{X_r}},~~\forall \phi \in {T_{{v_k^{(n)}}}} := \{ \phi  \in X_r, { \langle {{v_k^{(n)}},\phi}  \rangle _{{L^2}}} = 0\}.
$$
Indeed, taking $\widetilde\phi (x)=e^{ \frac{  {-N \theta_k^{(n)} }  }{2}  } \phi( e^{-\theta_k^{(n)}} x)$ for $\phi \in {T_{{v_k^{(n)}}}}$, we have
$$\begin{aligned}
 \langle {{I^{'}}({v_k^{(n)}}),\phi }  \rangle
=& \int_{{\R^N}} {\nabla {v_k^{(n)}} \nabla \phi } dx + \int_{{\R^N}} { | {\nabla {v_k^{(n)}} }  |}^{q-2} {\nabla {v_k^{(n)}} \nabla \phi }dx-\int_{{\R^N}} {{{ | {v_k^{(n)}} |}^{p-2}} {v_k^{(n)}} \phi }dx  \nonumber \\
=& {e^{ 2{\theta _k^{(n)}} }}\int_{{\R^N}} {\nabla {u_k^{(n)}} \nabla \widetilde \phi } dx + {e^{ q(1+\delta_q){\theta _k^{(n)}} }}\int_{{\R^N}} { | {\nabla {u_k^{(n)}} } |}^{q-2} {\nabla {u_k^{(n)}} \nabla \widetilde \phi }dx \nonumber \\
 &-{e^{ p \delta_p {\theta _k^{(n)}} }}\int_{{\R^N}} {{{ | {u_k^{(n)}} |}^{p-2}} {u_k^{(n)}} \widetilde \phi }dx  \nonumber \\
=&  \langle {{{\widetilde I}^{'}}({u_k^{(n)}},{\theta _k^{(n)}}),(\widetilde \phi ,0)}  \rangle.
\end{aligned}$$
As $\int_{{\R^N}} {{u_k^{(n)}}\widetilde \phi }dx=\int_{{\R^N}} {{v_k^{(n)}}\phi }dx$, we see that
$$
\phi  \in {T_{{v_k^{(n)}}}} \Leftrightarrow (\widetilde \phi ,0) \in {\widetilde T_{({u_k^{(n)}},{\theta _k^{(n)}})}}.
$$
Hence, utilizing $(b)$, we have
\[
 | {{\theta _k^{(n)}}}  | =  | {{\theta _k^{(n)}} - 0} | \le \mathop {\min }\limits_{0 \le t \le 1,u \in {S_{r,c}}  \cap {V_n}} { \| {({u_k^{(n)}},{\theta _k^{(n)}}) - ({\gamma_k}(t,u),0)}  \|_E} \le \frac{1}{{\sqrt k }}.
\]
The above estimate on $\theta _k^{(n)}$ indicates that
\[
 \| {(\widetilde \phi, 0)}  \|_E = { \| {\widetilde \phi }  \|}_{X_r}= {e^{{-\theta _k^{(n)}} }}  {{\| {\nabla \phi } \|}_2}+{e^{ -(1+\delta_q){\theta _k^{(n)}} }}{{\| {\nabla \phi } \|}_q}+{{\| { \phi } \|}_2} \le 2{\left\| \phi  \right\|}_{X_r}.
\]
Thus, we have
\[
| { \langle {{I^{'}}({v_k^{(n)}}),\phi } \rangle } | = \big \langle {{{\widetilde I}^{'}}({u_k^{(n)}},{\theta _k^{(n)}}),(\widetilde \phi,0)}  \big \rangle  \le \frac{2}{{\sqrt k }} \| {(\widetilde\phi,0)}  \|_E \le \frac{4}{{\sqrt k }}{ \| \phi   \|}_{X_r}.
\]
Consequently, we get
\[
\| {({\left. {{I}} \right|}_{{S_{r,c}} })^{'}({v_k^{(n)}})} \| = \mathop {\sup }\limits_{\phi  \in {T_{{v_k^{(n)}}},~\| \phi \| \le 1} }    | {\langle {{I^{'}}({v_k^{(n)}}),\phi } \rangle } | \le \frac{4}{{\sqrt k }} \to 0 \text{ as }k \to \infty .
\]
\end{proof}

To prove the compactness of the Palais-Smale sequence $\{ {v_k^{(n)}}\}$ obtained in Proposition \ref{5.7-ps}, we need the following Lemma.

\begin{lemma} \label{5.8-pq}
Consider the functional $I\!\in\!{C^1(X_r,\R)}$ and let $\left\{ {{u_k}} \right\} \!\subset\! S_{r,c}$ be a bounded sequence in $X_r$, then there exists $\lambda_k=\frac{\langle {I^{'}}\left( {{u_k}} \right),{u_k}\rangle}{c^2}$ such that
\[
(\left. {{I}} \right|_{ S_{r,c} })^{'}\left( {{u_k}} \right)\to 0 {\text{ in }}  X_r^{ - 1}  \Longleftrightarrow {I^{'}}\left( {{u_k}} \right) -\lambda_k {u_k} \to 0 {\text{ in }} X_r^{ - 1} \text{ as }k\rightarrow +\infty.
\]
\end{lemma}

\begin{proof}
The proof is similar to that of \cite[Lemma 3]{HbPl}.
\end{proof}

\begin{proposition}\label{5.9-ps}
Assume \eqref{h1.3a}, \eqref{h1.3b}. Let $\{ {v_k^{(n)}}\}\subset {S_{r,c}}$ be the Palais-Smale sequence obtained in Proposition \ref{5.7-ps}, then there exists ${v^{(n)}} \in X_r$ and $\{\lambda _k^{(n)}\}\subset \R$ such that, up to a subsequence, \\
$(i)~~ {v_k^{(n)}} \rightharpoonup {v^{(n)}} \not \equiv 0\text{ in }X_r$ as $k \to +\infty$;\\
$(ii)~~ \lambda_k^{(n)} \to \lambda^{(n)}< 0 \text{ for some } \lambda^{(n)}\in\R$ as $k \to +\infty$;\\
$(iii)~~ -\Delta {v_k^{(n)}}- \Delta_q {v_k^{(n)}} - {\lambda _k^{(n)}}{v_k^{(n)}} - {| {{v_k^{(n)}}} |^{p - 2}}{v_k^{(n)}} \to 0\text{ in } X_r^{-1}$ as $k \to +\infty$;\\
$(iv)~~ - \Delta {v^{(n)}}- \Delta_q {v^{(n)}}-{\lambda^{(n)}}{v^{(n)}} - {| {{v^{(n)}}} |^{p-2}}{v^{(n)}}=0\text{ in } X_r^{-1}$;\\
$(v)~~{v_k^{(n)}} \to {v^{(n)}} \text{ in } X_r\text{ as }k \to +\infty$.
\end{proposition}

\begin{proof}
By using Lemma \ref{5.7-ps}, we have that $P(v_k^{(n)})\to 0$ as $k\to\infty$, so that, for $k$ sufficiently small, we get
$$\biggl(\frac{1}{2}-\frac{1}{p\delta_{p}}\biggr) \|\nabla v_k^{(n)}\|_{2}^2+\biggl(\frac{1}{q}-\frac{1+\delta_{q}}{p\delta_{p}}\biggr)\|\nabla v_k^{(n)}\|_{q}^q=I(v_k^{(n)})\leq \sigma_n(c)+1.$$
Therefore $\{ {v_k^{(n)}}\} $ is bounded in $X_r$ and there exists ${v^{(n)}} \in X_r$ such that
$$
 { v_k^{(n)} } \rightharpoonup {v^{(n)}}\text{ in } X_r,~~~~
 { v_k^{(n)} } \to {v^{(n)}}\text{ in }{L^p_{loc}}({\R^N}),~~~~
 { v_k^{(n)} } \to {v^{(n)}}\text{ a.e on }{\R^N},
$$
up to a subsequence. We claim that $v^{(n)}\not \equiv 0$. Otherwise, ${\| { {v_k^{(n)} } } \|}_p\to0$, from $P({v_k^{(n)}})=o_k(1)$, we have
$$
{{\| {\nabla {v_k^{(n)}} } \|}_2^2}+  (1+\delta_q){{\| {\nabla {v_k^{(n)}} } \|}_q^q}=\delta_p  {\| { {v_k^{(n)} } } \|}_p^p+o_k(1)=o_k(1),
$$
which leads to $I({v_k^{(n)}})=o_k(1)$. However, this contradicts to $\sigma_n(c)\ge {\beta_n}>1$. Thus, $(i)$ is true. From Lemma \ref{5.8-pq}, there exists ${\lambda_k^{(n)}}\in\R$ such that
\[
{I^{'}} ( {{v_k^{(n)}}}) - {\lambda_k^{(n)}} {v_k^{(n)}} \to 0 {\text{ in }}{X_r^{ - 1}}
\]
as $k\rightarrow\infty$, where \begin{equation}\label {lam5.13}
\begin{array}{rl}
{\lambda_k^{(n)}}=\dfrac{\langle {I^{'}}( {{v_k^{(n)} }} ),{v_k^{(n)} }\rangle}{c^2}=\dfrac{  {{\| {\nabla {v_k^{(n)}} } \|}_2^2}+  {{\| {\nabla {v_k^{(n)}} } \|}_q^q}- {\| { {v_k^{(n)} } } \|}_p^p   }{ c^2  }.
\end{array}
\end{equation}
That is to say, for any $\phi\in {X_r}$, we have
\begin{align}\label{5.14-nk}
\Big\langle {I^{'}}( {{v_k^{(n)} }} ) - {\lambda_k^{(n)}}  {v_k^{(n)}},\phi \Big \rangle
=&\int_{{\R^N}} {\nabla {v_k^{(n)}} \nabla \phi } dx + \int_{{\R^N}} { | {\nabla {v_k^{(n)}} }  |}^{q-2} {\nabla {v_k^{(n)}} \nabla \phi }dx  \nonumber \\
&-\int_{{\R^N}} {{{ | {v_k^{(n)}} |}^{p-2}} {v_k^{(n)}} \phi }dx -{\lambda _k^{(n)}}\dis \int_{{\R^N}} {{v_k^{(n)}}\phi }dx=o_k(1)\|\phi\|_{X_r}.
\end{align}
Hence $(iii)$ holds. The boundedness of $\{ {v_k^{(n)}}\} $ in $X_r$ and $P({v_k^{(n)}}) =o_k(1)$ imply that $\{ \lambda_k^{(n)} \}$ is bounded. Then, there exists $\lambda^{(n)} \in \R$ such that $\lambda_k^{(n)} \to \lambda^{(n)}$ as $k \to +\infty$ up to a subsequence. In addition, taking into account \eqref{lam5.13}, then $P({v_k^{(n)}})=o_k(1)$ indicates that
\begin{align}  \label{5.15lgrange}
\lambda^{(n)} &= \mathop {\lim }\limits_{k \to \infty }  \lambda_k^{(n)}= \mathop {\lim }\limits_{k \to \infty } \dis\frac{1}{c^2}[ {{\| {\nabla {v_k^{(n)}} } \|}_2^2}+  {{\| {\nabla {v_k^{(n)}} } \|}_q^q}- {\| { {v_k^{(n)} } } \|}_p^p  ] \nonumber \\
&= \mathop {\lim }\limits_{k \to \infty } \dis\frac{1}{c^2}\biggl[\biggl(1-\frac{1}{\delta_p}\biggr) {{\| {\nabla {v_k^{(n)}} } \|}_2^2}+\biggl(1-\frac{1+\delta_q}{\delta_p}\biggr){{\| {\nabla {v_k^{(n)}} } \|}_q^q} \biggr ].
\end{align}
Notice that $1-\frac{1}{\delta_{p}}<0$ and $1-\frac{ 1+\delta_{q} }{\delta_{p}}<0$, we have $\lambda^{(n)}\leq0$. It must be that $\lambda^{(n)}<0$. Otherwise, if $\lambda^{(n)}=0$, then \eqref{5.15lgrange} gives
$$
\mathop {\lim }\limits_{k \to \infty } {{\| {\nabla {v_k^{(n)}} } \|}_2}=\mathop {\lim }\limits_{k \to \infty }{{\| {\nabla {v_k^{(n)}} } \|}_q}=0,
$$
but this contradicts to $\sigma_n(c)>1$. Thus, $(ii)$ is true. From $(ii)-(iii)$, then $v^{(n)}$ satisfies
\begin{align}\label {5.16-n}
\int_{{\R^N}}\Big [ {\nabla {v^{(n)}} \nabla \phi }   +   { | {\nabla {v^{(n)}} }  |}^{q-2} {\nabla {v^{(n)}} \nabla \phi } -  {{{ | {v^{(n)}} |}^{p-2}} {v^{(n)}} \phi } -{\lambda^{(n)}}   {{v^{(n)}}\phi } \Big]dx=0,~~~~\forall \phi\in {X_r},
\end{align}
That is equivalent to $(iv)$. Testing \eqref{5.14-nk} and \eqref{5.16-n} with $\phi=v_k^{(n)}-v^{(n)}$ and using ${v_k^{(n)}} \to {v^{(n)}}\text{ in }{L^p}({\R^N})$, we have
\begin{equation}\label{fin}
\begin{aligned}
\int_{{\R^N}}& {{{ | {\nabla ({v_k^{(n)}} \!-\!  {v^{(n)}})} |}^2}} \!+ \! {\big({{ | {\nabla {v_k^{(n)}}}  |}^{q-2}}{\nabla v_k^{(n)}} \!- \!{{ | \nabla {{v^{(n)}}}  |}^{q - 2}}{\nabla v^{(n)}}\big )\nabla ( {v_k^{(n)}} \!- \! {v^{(n)}})} \!- \!{\lambda ^{(n)}} {{{ | {{v_k^{(n)} } - {v^{(n)}}}  |}^2}} dx \\
&
= \int_{{\R^N}} {({{ | {{v_k^{(n)}}}  |}^{p - 2}}{v_k^{(n)}} - {{ | {{v^{(n)}}} |}^{p - 2}}{v^{(n)}})({v_k^{(n)}} - {v^{(n)}})}  + o_k(1)=o_k(1).
\end{aligned}
\end{equation}
For any $\xi,\eta\in\R^N $ with $|\xi|+|\eta|>0$, we learn from \cite[Lemma 2.1]{LdLi} that the following inequalities hold
\begin{align*}
({ | \xi | }^{q - 2}{\xi} -  { | \eta |}^{q - 2} { \eta })( \xi-\eta)&\geq C|\xi-\eta|^q,~~~~q>2   \\
({ | \xi | }^{q - 2}{\xi} -  { | \eta |}^{q - 2} { \eta })( \xi-\eta)&\geq C(|\xi|+|\eta|)^{q-2}|\xi-\eta|^2,~~~~1<q\leq2
\end{align*}
for some constant $C>0$. Hence, if $q>2$, from \eqref{fin}, we have
$$\begin{aligned}
 &\int_{{\R^N}} \Big [{{{ | {\nabla ({v_k^{(n)}} \!-\!  {v^{(n)}})} |}^2}} \!+ \!C {{{ | {\nabla ({v_k^{(n)}} \!-\!  {v^{(n)}})} |}^q}}  \!- \!{\lambda ^{(n)}} {{{ | {{v_k^{(n)} } \!-\! {v^{(n)}}}  |}^2}} \Big] dx =o_k(1).
\end{aligned}$$
If $1<q<2$, we deduce that
$$C^{ \frac{q}{2} } |\xi-\eta|^q \leq  \big |({ | \xi | }^{q - 2}{\xi} -  { | \eta |}^{q - 2} { \eta })( \xi-\eta) \big |^{ \frac{q}{2} }(|\xi|+|\eta|)^{\frac{q(2-q)}{2}}.$$
In view of this inequality, we use the H\"{o}lder's inequality to obtain
$$\begin{aligned}
& C^{ \frac{q}{2} } \int_{{\R^N}}   {{{ | {\nabla ({v_k^{(n)}} \!-\!  {v^{(n)}})} |}^q}}   dx  \leq
   \int_{{\R^N}} \big | B(v_k^{(n)},v^{(n)}) \big |^{ \frac{q}{2} }  (|\nabla {v_k^{(n)}}|+|\nabla {v^{(n)}}|)^{\frac{q(2-q)}{2}} dx \\ \nonumber
 \leq&
   \Big \{ \int_{{\R^N}} \big | B(v_k^{(n)},v^{(n)}) \big |  dx \Big \}^{ \frac{q}{2} } \Big \{ \int_{{\R^N}}  (|\nabla {v_k^{(n)}}|+|\nabla {v^{(n)}}|)^{q} dx \Big \}^{ \frac{2-q}{2} } \\ \nonumber
 \leq&
   \Big \{ \int_{{\R^N}} \big | B(v_k^{(n)},v^{(n)})  \big |  dx \Big \}^{ \frac{q}{2} }   \Big \{ \int_{{\R^N}}  2^{q-1}(|\nabla {v_k^{(n)}}|^{q}+|\nabla {v^{(n)}}|^{q}) dx \Big \}^{ \frac{2-q}{2} },
\end{aligned}$$
where $B(v_k^{(n)},v^{(n)})\!=\!({ | \nabla {v_k^{(n)}} | }^{q - 2}{ \nabla {v_k^{(n)}} } -  { | \nabla {v^{(n)}} |}^{q - 2} { \nabla {v^{(n)}} })  \nabla ({v_k^{(n)}}-  {v^{(n)}})$. As $\{v_k^{(n)}\} $ is bounded in $X_r $, we have
$$
 \int_{{\R^N}}   ({ | \nabla {v_k^{(n)}} | }^{q - 2}{ \nabla {v_k^{(n)}} } -  { | \nabla {v^{(n)}} |}^{q - 2} { \nabla {v^{(n)}} })  \nabla ({v_k^{(n)}}-  {v^{(n)}})   dx \geq C \Big \{  \int_{{\R^N}}   {{{ | {\nabla ({v_k^{(n)}} \!-\!  {v^{(n)}})} |}^q}}   dx  \Big \}^{ \frac{2}{q} }.
$$
Consequently, if $1<q<2$, from \eqref{fin}, we have
$$
   \| {\nabla ({v_k^{(n)}} \!-\!  {v^{(n)}})} \|_2^2  \!+ \!C  \| {\nabla ({v_k^{(n)}} \!-\!  {v^{(n)}})} \|_q^2  \!- \!{\lambda ^{(n)}}  \| {{v_k^{(n)} } \!-\! {v^{(n)}}}  \|_2^2=o_k(1).
$$
Recalling that $\lambda^{(n)}<0$, so we obtain $v_k^{(n)} \to v^{(n)}$ in $X_r $ as $k\to \infty$. Thus $(v)$ is proved.
\end{proof}

\begin{proof}[\textbf {Proof of Theorem \ref{th1.4}.} ]
For any fixed $c>0$ and $n\geq n_0$, Proposition \ref{5.9-ps} implies that $\eqref{eq1.1}$ possesses a sequence of weak solutions $\{({v^{(n)}},{\lambda^{(n)}})\}\subseteq X_r \times {\R^- }$ with $\left\| {{v^{(n)}}} \right\|_2=c$. Since $P({v^{(n)}})= 0$, from Lemma \ref{6.1.sig-beta} we have
\begin{align*}
\biggl(\frac{1}{2}-\frac{1}{p\delta_{p}}\biggr) \|\nabla v^{(n)}\|_{2}^2+\biggl(\frac{1}{q}-\frac{1+\delta_{q}}{p\delta_{p}}\biggr)\|\nabla v^{(n)}\|_{q}^q=I(v^{(n)})= \sigma_n(c)\ge {\beta_n} \to +\infty
\end{align*}
as $n \to +\infty$. Consequently, we have $I(v^{(n)})\to +\infty$ and $\|v^{(n)}\|_{X_r}\to +\infty$.
\end{proof}

\vspace{.55cm}
\noindent{\bf Acknowledgements:}
Tao Yang was supported by National Natural Science Foundation of China (Grant No. 12201564). Laura Baldelli is a member of the {\em Gruppo Nazionale per l'Analisi Ma\-te\-ma\-ti\-ca, la Probabilit\`a e le loro Applicazioni} (GNAMPA) of the {\em Istituto Nazionale di Alta Matematica} (INdAM).
Laura Baldelli was partially supported by National Science Centre, Poland (Grant No. 2020/37/B/ST1/02742).

\end{document}